\documentclass[12pt]{amsart}
\usepackage{amsmath,amscd,amssymb,amsfonts,latexsym}
\usepackage[all]{xy}

\usepackage{geometry}
\allowdisplaybreaks

\newcommand{\bC}{{\mathbb C}}

\newcommand{\bQ}{{\mathbb Q}}
\newcommand{\bZ}{{\mathbb Z}}
\newcommand{\bV}{{\mathbb V}}

\newcommand{\cD}{{\mathcal D}}

\newcommand{\cF}{{\mathcal F}}
\newcommand{\cG}{{\mathcal G}}
\newcommand{\cH}{{\mathcal H}}

\newcommand{\cL}{{\mathcal L}}
\newcommand{\cM}{{\mathcal M}}

\newcommand{\cO}{{\mathcal O}}

\newcommand{\cS}{{\mathcal S}}

\newcommand{\So}{\,\overline{\!S}}

\newcommand{\Gr}{\text{\rm Gr}}
\newcommand{\DR}{\hbox{\rm DR}}

\newcommand{\h}{\hbox}
\newcommand{\q}{\quad}

\newcommand{\os}{$O_{\sigma}$}
\newcommand{\vs}{$V_{\sigma}$}

\newcommand{\tyh}{$\widehat{T}_{y*}$}

\newcommand{\La}{\Lambda}

\def\sig{\sigma}
\def\Sig{\Sigma}

\newcommand{\MHM}{{\rm MHM}}
\newcommand{\MHS}{{\rm MHS}}

\newcommand{\VMHS}{{\rm VMHS}}

\newcommand{\rat}{{\rm rat}}
\newcommand{\var}{{\rm var}}
\newcommand{\Perv}{{\rm Perv}}
\newcommand{\im}{{\rm Image}}

\newcommand{\Xs}{{X^{(n)}}}

\theoremstyle{plain}
\newtheorem{thm}{Theorem}[section]
\newtheorem{cor}[thm]{Corollary}

\newtheorem{lem}[thm]{Lemma}
\newtheorem{prop}[thm]{Proposition}

\theoremstyle{definition}
\newtheorem{df}[thm]{Definition}
\newtheorem{rem}[thm]{Remark}

\newtheorem{example}[thm]{Example}

\def\be{\begin{equation}}
\def\ee{\end{equation}}

\def\bt{\begin{thm}}
\def\et{\end{thm}}

\def\bc{\begin{cor}}
\def\ec{\end{cor}}

\def\br{\begin{rem}}
\def\er{\end{rem}}

\def\bp{\begin{prop}}
\def\ep{\end{prop}}

\def\bl{\begin{lem}}
\def\el{\end{lem}}

\def\bn{\begin{enumerate}}
\def\en{\end{enumerate}}

\def\bex{\begin{example}}
\def\eex{\end{example}}

\def\bd{\begin{df}}
\def\ed{\end{df}}

\numberwithin{equation}{section}

\title{Characteristic classes of mixed Hodge modules\\ and applications}
\author{Lauren\c{t}iu Maxim}
\address{L. Maxim: Department of Mathematics, University of Wisconsin, 480 Lincoln Drive, Madison, WI 53706, USA}
\email {maxim@math.wisc.edu}
\author[J. Sch\"urmann ]{J\"org Sch\"urmann}
\address{J.  Sch\"urmann : Mathematische Institut,
          Universit\"at M\"unster,
          Einsteinstr. 62, 48149 M\"unster,
          Germany.}
\email {jschuerm@math.uni-muenster.de}

\date{\today}



\begin{document}

\begin{abstract}
These notes are an extended version of the authors' lectures at the 2013 CMI Workshop ``Mixed Hodge Modules and Their Applications". We give an overview, with an emphasis on applications, of recent developments on the interaction between characteristic class theories for singular spaces and Saito's theory of mixed Hodge modules in the complex algebraic context.
\end{abstract}

\maketitle

\tableofcontents

\section{Introduction}\label{intro}

We give an overview of recent developments on
the interaction between characteristic class theories for singular spaces and Saito's theory of mixed Hodge modules in the complex algebraic context, updating the existing survey \cite{Sc}. The emphasis here is on applications.\\
      
There are two versions of characteristic classes associated to mixed Hodge modules, cf. \cite{BSY,Sc}. The
$K$-theoretical classes, called {\it Hodge-Chern classes}, capture information about the graded pieces of the filtered de Rham complex associated to the filtered $\cD$-module underlying a
mixed Hodge module. Application of the Todd class
transformation $td_*$ of Baum-Fulton-MacPherson then gives classes in (Borel-Moore) homology $H_*:=H_{2*}^{BM}$, which we call {\it un-normalized Hirzebruch classes}. Both versions are defined by natural transformations $$\DR_y:K_0(\MHM(X)) \to K_0(X) \otimes \bZ[y^{\pm 1}]$$ and resp. $$T_{y*}:=td_* \circ \DR_y :K_0(\MHM(X)) \to H_*(X)[y,y^{-1}]$$ on the Grothendieck group of mixed Hodge modules on a variety $X$. For the {\it normalized Hirzebruch classes} $\widehat{T}_{y*} $ we renormalize the classes $T_{y*}$ by multiplying by $(1+y)^{-i}$ on $H_i(-)$.

These transformations are
functorial for proper pushdown and external products, and satisfy
some other properties which one would expect for a theory of
characteristic classes for singular spaces.
For ``good"  variations of mixed Hodge structures (and their extensions along a normal crossing divisor) the corresponding Hodge-Chern and resp. Hirzebruch classes have an
explicit description in terms of (logarithmic) de Rham
complexes. On a point space, these classes coincide with
the {\it Hodge polynomial} $\chi_y$ of a mixed Hodge structure. From this point of view, the Hirzebruch characteristic classes of mixed Hodge modules can be regarded as higher homology class versions of the Hodge polynomials defined in terms of the Hodge filtration on the cohomology of the given mixed Hodge module.
This will be explained in Section \ref{HCH}, after some relevant background about mixed Hodge modules is presented in Section \ref{mixedHodge}.
As a first application, these characteristic classes are computed in the case of toric varieties, see Section \ref{Toric}.\\

In Section \ref{HM}, we discuss the specialization property of Hodge-Chern and resp. Hirzebruch classes, which can be viewed as a Hodge-theoric counterpart of the Verdier specialization result for the MacPherson Chern class. 
In more detail, for  a globally defined hypersurface $X=\{f=0\}$ given by the zero-fiber of a complex algebraic function $f:M \to \bC$ defined on an algebraic manifold $M$,  it is proved in \cite{Sb} that 
the Hodge-Chern and resp. normalized Hirzebruch class transformations commute with specialization defined in terms of the nearby cycles $\psi^H_f$ with 
$\rat(\psi^H_f)={^p}\psi_f$, in the sense that 
$$(1+y)\cdot \DR_y \circ (\psi^H_f[1])=i^!\circ \DR_y$$ and
$$\widehat{T}_{y*} \circ (\psi^H_f[1])=i^!\circ \widehat{T}_{y*},$$ for $i:X \to M$ the inclusion map and $i^!$ the corresponding Gysin homomorphism.
The proof of this result, discussed in Section \ref{specv}, uses the algebraic theory of nearby and vanishing cycles given 
by the $V$-filtration of Malgrange-Kashiwara (as recalled in Section \ref{v}) in the $\cD$-module context, together with a specialization result about the filtered de Rham complex of the filtered $\cD$-module underlying 
a mixed Hodge module.

As an application, we give in Section \ref{applic} a description of the difference between
the corresponding virtual and functorial Hirzebruch characteristic classes of complex
hypersurfaces in terms of vanishing cycles related to the singularities of
the hypersurface. See also \cite{MSS2} for the case of projective hypersurfaces
(and even projective complete intersections).\\


In Section \ref{equivariant}, we discuss equivariant characteristic class theories for singular varieties.
More precisely, following  \cite{CMSS2}, for varieties with {\em finite} algebraic group actions we define for a ``weak" equivariant complex of mixed Hodge modules
localized {\it Atiyah-Singer classes} in the (Borel-Moore) homology of the fixed point sets.
These are equivariant versions of the above Hirzebruch classes, and enjoy similar functorial properties. They can be used for computing explicitly 
the Hirzebruch classes of global quotients, e.g., of symmetric products of quasi-projective varieties.

In Section \ref{symmhm}, we recall the definition of symmetric products of mixed Hodge module complexes. 
For a complex quasi-projective variety $X$, let $X^{(n)}:=X^n/\Sigma_n$ denote its $n$-th {\it symmetric product}, with $\pi_n:X^n \to X^{(n)}$ the natural 
projection map.
In \cite{MSS}, we define an action of the symmetric group $\Sigma_n$ on the  $n$-fold external self-product $\boxtimes^n\cM$ of an arbitrary 
bounded complex of mixed Hodge modules $\cM \in D^b\MHM(X)$. By construction, this action is compatible with the natural action on the underlying 
$\bQ$-complexes. There are, however,  certain technical difficulties associated with this construction, since the difference in the $t$-structures
of the underlying $\cD$-modules and $\bQ$-complexes gives certain differences of signs. In \cite{MSS} we solve this problem by showing that there is a sign 
cancellation. 
For a complex of mixed Hodge modules $\cM \in D^b\MHM(X)$, we can therefore define its {\it $n$-th symmetric product} 
(in a way compatible with $\rat$) by:
$${\cM^{(n)}}:=({\pi_n}_*\cM^{\boxtimes{n}})^{\Sigma_n} 
\in D^b\MHM(X^{(n)}).$$

In Section \ref{applsym} we discuss a generating series formula for the un-normalized Hirzebruch classes of symmetric products of a mixed Hodge module complex, which 
generalizes 
to this singular setting many of the generating series formulae in the literature. The key technical point  is a localization formula for the 
Atiyah--Singer classes, which relies on understanding how Saito's functors $Gr^F_p\DR$ behave with respect to external products.\\

\textbf{Acknowledgments.}  We thank J.-P. Brasselet, S. Cappell, A. Libgober, M. Saito, J. Shaneson, S. Yokura for many useful discussions and for collaborating with us on various parts of the research described in this note.


\section{Calculus of Mixed Hodge Modules}\label{mixedHodge}

For the sake of completeness and coherence of exposition, in this section we give a brief overview of the theory of mixed Hodge modules and their Grothendieck calculus.

\subsection{Mixed Hodge Modules}

To any complex algebraic variety $X$, Saito associated a
category $\MHM(X)$ of {\em algebraic mixed Hodge modules} on $X$ (cf.
\cite{Sa0,Sa1}). If $X$ is smooth, an object of this category
consists of an algebraic (regular) holonomic $\cD$-module $(\cM,F)$ with a good filtration $F$,
together with a perverse sheaf $K$ of rational vector spaces, both endowed a finite
increasing filtration $W$ such that the isomorphism 
$$\alpha: \DR(\cM)^{an}\simeq K\otimes_{\bQ_X} \bC_X $$
is compatible with $W$
under the Riemann-Hilbert correspondence (with $\alpha$ a chosen isomorphism). Here we use {\it left} $\cD$-modules. The sheaf $\cD_X$ of algebraic differential operators on $X$ has the increasing filtration $F$ with $F_i\cD_X$ given by the differential operators of degree
$\leq i$ ($i\in \bZ$). A {\em good} filtration $F$ of the algebraic holonomic $\cD$-module $\cM$ is then given by a bounded from below, increasing, and exhaustive filtration $F_p\cM$ by {coherent} algebraic $\cO_X$-modules such that
\begin{equation}\label{D-filt}
F_i\cD_X\left( F_p\cM\right) \subset F_{p+i}\cM 
\end{equation}
for all $i,p$, and this is an equality for $i$ big enough.

In general, for a singular variety $X$ one works
with suitable local embeddings into manifolds and corresponding
filtered $\cD$-modules supported on $X$. In addition, these objects
are required to satisfy a long list of  properties (which are not needed here).

The {\em forgetful functor} $\rat$ is defined as
$$\rat:  \MHM(X)\to \Perv(\bQ_X) \ , \ \ \ \ \left((\cM,F),K,W\right)\mapsto K\:,$$
with $\Perv(\bQ_X)$ the abelian category of perverse sheaves on $X$.

For the following result, see \cite{Sa1}[Thm.0.1 and Sec.4] for more details:
\bt [M. Saito]\label{MHM1}
The category $\MHM(X)$ is abelian, and the forgetful functor $\rat:  \MHM(X)\to \Perv(\bQ_X)$
is exact and {faithful}. It extends to a functor 
$$\rat: D^b\MHM(X) \to D^b_c(\bQ_X)$$ 
to the derived category of complexes of $\bQ$-sheaves with algebraically constructible
cohomology. There are functors 
$$f_*, \;f_! ,\; f^*,\; f^!,\; \otimes,\; \boxtimes  $$ 
on $D^b\MHM(X)$
which are
``lifts" via $\rat$ of the similar (derived) functors defined on $D^b_c(\bQ_X)$, and with $(f^*,f_*)$ and $(f_!,f^!)$ pairs of {adjoint} functors.
One has a natural map $f_!\to f_*$, which is an isomorphism for $f$ proper.
\et
The usual truncation $\tau_{\leq}$ on $D^b\MHM(X)$
corresponds to the {\em perverse truncation} ${^p\tau}_{\leq}$ on
$D^b_c(X)$, so that
$$\rat\circ H= \;^p\cH\circ \rat\:,$$
where $H$ stands for the cohomological functor in $D^b\MHM(X)$ and $\;^p\cH$ denotes the perverse cohomology
(with respect to the middle perversity).

\begin{example}\label{smooth}
Let $X$ be a complex algebraic manifold of pure complex dimension $n$,
with $\bV:=(L,F,W)$ a {\it good} (i.e., admissible, with quasi-unipotent monodromy at infinity) variation of mixed Hodge structures on $X$.
Then $\cL:=L \otimes_{\bQ_X} \cO_X$ with its integrable connection $\nabla$ is a holonomic (left) $\cD$-module
with $\alpha: \DR(\cL)^{an}\simeq L[n]$, where we use the shifted de Rham complex
$$\begin{CD}
\DR(\cL):=[\cL @> \nabla >> \cdots @> \nabla >> \cL\otimes_{\cO_X} \Omega^n_X]
\end{CD}$$ 
with $\cL$ in degree $-n$, so that $DR(\cL)^{an}\simeq L[n]$ is a perverse sheaf on $X$.
The filtration $F$ on $\cL$ induces by Griffiths' transversality a good filtration
$F_p(\cL):=F^{-p}\cL$ on $\cL$ as a filtered  $\cD$-module. 
Moreover,  $\alpha$ is compatible with the
induced filtration $W$ defined by 
$$W^i(L[n]):=W^{i-n}L[n] \quad \text{ and} \quad  
W^i(\cL):=(W^{i-n}L)\otimes_{\bQ_X} \cO_X \:.$$
 This data defines a mixed Hodge module
$\bV^H[n]$ on $X$, with $\rat(\bV^H[n])\simeq L[n]$. Hence $\rat(\bV^H[n])[-n]$ is a local system on $X$.  

In what follows, we will often use the same symbol $\bV$ to denote both the variation and the corresponding (shifted) mixed Hodge module.
\end{example}

\bd A mixed Hodge module $\cM$ on a pure $n$-dimensional complex algebraic manifold
$X$ is called {\em smooth} if $rat(\cM)[-n]$ is a local system on $X$.\ed

We have the following result (see \cite{Sa1}[Thm.3.27, Rem. on p.313]):
\bt[M. Saito]\label{MHM2}
Let $X$ be a pure $n$-dimensional complex algebraic manifold.
Associating to a good variation of mixed Hodge structures $\bV=(L,F,W)$ on $X$
the mixed Hodge module $\bV^H[n]$ as in Example (\ref{smooth}) defines an
equivalence of categories
$$\MHM(X)_{sm}\simeq \VMHS^g(X)$$
between the categories of smooth mixed Hodge modules $\MHM(X)_{sm}$ and 
good variations of mixed Hodge structures on $X$.
This equivalence commutes with external product $\boxtimes$. 
For $X=pt$ a point, one obtains in particular an equivalence
\be\label{1} \MHM(pt)\simeq \MHS^p \ee
between mixed Hodge modules on a point space and the abelian category of graded polarizable mixed Hodge structures.
\et

By the identification in (\ref{1}), there exists a unique Tate object  $\bQ^H(k)
\in \MHM(pt)$ such that $\rat(\bQ^H(k))=\bQ(k)$, with $\bQ(k)$ of type $(-k,-k)$. For a complex variety $X$ with constant map $k: X\to pt$, define 
$$\bQ_X^H(k):=k^*\bQ^H(k) \in D^b\MHM(X), $$ with $$\rat(\bQ_X^H(k))=\bQ_X(k).$$ 
So tensoring with $ \bQ_X^H(k)$ defines the {\it Tate twist operation} $\cdot (k)$  on mixed Hodge modules. 
To simplify the notations, we let $\bQ_X^H:=\bQ_X^H(0)$.
If $X$ is \emph{smooth} of complex dimension $n$ then
$\bQ_X[n]$ is perverse on $X$, and $\bQ_X^H[n]\in \MHM(X)$ is a single
mixed Hodge module, explicitly described by
$$\bQ_X^H[n]=((\cO_X, F), \bQ_X[n], W), $$ with  $gr^F_i=0=gr^W_{i+n}$ for all $ i \neq 0$.\\

Let us also mention here the following result about the existence of a mixed Hodge structure
on the cohomology (with compact support) $H^i_{({c})}(X;\cM)$ for $\cM \in D^b\MHM(X)$.

\bc\label{MHM3}
Let $X$ be a complex algebraic variety with constant map $k: X\to pt$.
Then the cohomology (with compact support) $H^i_{({c})}(X;\cM)$ of $\cM \in D^b\MHM(X)$ gets an induced graded polarizable mixed Hodge structure by:
$$H^i_{({c})}(X,\cM)=H^i(k_{*(!)}\cM)\in \MHM(pt)\simeq \MHS^p\:.$$
In particular:
\begin{enumerate}
\item The rational cohomology (with compact support) $H^i_{({c})}(X;\bQ)$ 
of $X$ gets an induced graded polarizable mixed Hodge structure by:
$$H^i(X;\bQ)= \rat(H^i(k_*k^*\bQ^H)) \quad \text{and} \quad 
H^i_{{c}}(X;\bQ)= \rat(H^i(k_!k^*\bQ^H)) \:.$$
\item Let $\bV$ be a good variation of mixed Hodge structures on a
smooth (pure dimensional) complex manifold $X$. Then $H^i(X;\bV)$
gets a mixed Hodge structure by
$$H^i(X;\bV)\simeq \rat(H^i(k_*\bV^H))\:,$$
and, similarly, $H_c^i(X;\bV)$
gets a mixed Hodge structure by
$$H_c^i(X;\bV)\simeq \rat(H^i(k_!\bV^H))\:.$$
\end{enumerate}
\ec

\br 
 By a deep theorem of Saito (\cite{Sa3}[Thm.0.2,Cor.4.3]), the mixed Hodge structure on $H^i_{({c})}(X;\bQ)$ defined as above coincides with the classical mixed Hodge structure constructed by Deligne (\cite{De1, De3}).
 \er

We conclude this section with a short explanation of the {\it rigidity
property} for good variations of mixed Hodge structures. Assume $X$ is smooth, connected and of dimension $n$, with $\cM \in \MHM(X)$ a smooth mixed Hodge module, so
that the underlying local system $L:=\rat(\cM)[-n]$ has the property that the
restriction map $r:H^0(X;L) \to L_x$ is an isomorphism for all
$x \in X$. Then the (admissible) variation of mixed Hodge structures $\bV:=(L,F,W)$ 
is a {\it constant} variation since $r$ underlies the morphism of mixed
Hodge structures (induced by the adjunction $id \to i_*i^*$):
$$H^0(k_*\cM[-n]) \to H^0(k_*i_*i^*\cM[-n])$$
with $k:X \to pt$ the constant map, and $i :\{x\} \hookrightarrow X$
the inclusion of the point. This implies $$\cM[-n]=\bV^H\simeq
k^*i^*\bV^H=k^*\bV^H_x \in D^b\MHM(X).$$


\subsection{Grothendieck Groups of Algebraic Mixed Hodge Modules.}\label{Grot} 
In this section, we describe the functorial calculus of Grothendieck groups
of algebraic mixed Hodge modules. Let $X$ be a complex algebraic
variety. By associating to (the class of) a complex the alternating
sum of (the classes of) its cohomology objects, we obtain the
following identification
\begin{equation} K_0(D^b\MHM(X))=K_0(\MHM(X))
\end{equation}
between the corresponding Grothendieck groups.
In particular, if $X$ is a point, then
\begin{equation} K_0(D^b\MHM(pt))=K_0(\MHS^p),
\end{equation}
and the latter is a commutative ring with respect to the tensor
product, with unit $[\bQ^H]$.  
For any complex $\cM \in D^b\MHM(X)$, we have  the
identification
\begin{equation}\label{i1}
[\cM]=\sum_{i \in \bZ} (-1)^i [H^i(\cM)] \in
K_0(D^b\MHM(X)) \cong K_0(\MHM(X)).
\end{equation}
In particular, if for any $\cM \in \MHM(X)$ and $k \in \bZ$ we regard
$\cM[-k]$ as a complex concentrated in degree $k$, then
\begin{equation}\label{i2}
\left[ \cM[-k] \right]= (-1)^k [\cM] \in K_0(\MHM(X)).
\end{equation}
All functors $f_*$, $f_!$, $f^*$, $f^!$, $\otimes$, $\boxtimes$
induce corresponding functors on $K_0(\MHM(-))$. Moreover,
$K_0(\MHM(X))$ becomes a $K_0(\MHM(pt))$-module, with the
multiplication induced by the exact external product with a point space:
$$\boxtimes : \MHM(X) \times \MHM(pt) \to \MHM(X \times \{pt\}) \simeq
\MHM(X).$$ Also note that 
$$\cM \otimes \bQ^H_X \simeq \cM \boxtimes
\bQ^H_{pt} \simeq \cM$$ 
for all $\cM \in \MHM(X)$. Therefore,
$K_0(\MHM(X))$ is a unitary $K_0(\MHM(pt))$-module. The functors
$f_*$, $f_!$, $f^*$, $f^!$ commute with external products (and $f^*$
also commutes with the tensor product $\otimes$), so that the
induced maps at the level of Grothendieck groups $K_0(\MHM(-))$
are $K_0(\MHM(pt))$-linear. Moreover, by using the forgetful functor
$$\rat:K_0(\MHM(X)) \to K_0(D^b_c(\bQ_X)) \simeq K_0(\Perv(\bQ_X)),$$ all these transformations
lift the corresponding transformations from the (topological) level
of Grothendieck groups of constructible (or perverse) sheaves.

\br\label{generator}
The Grothendieck group $K_0(\MHM(X))$ is generated  by the classes $f_*[j_*\bV]$ (or, alternatively, by the classes $f_*[j_!\bV]$), with $f: Y\to X$ a proper morphism from a  complex algebraic manifold $Y$, $j: U\hookrightarrow Y$ the inclusion of a Zariski open and dense subset $U$, with complement $D$ a normal crossing divisor with smooth irreducible components, and $\bV$ a good variation of mixed (or pure) Hodge structures on $U$. 
This follows by induction from resolution of singularities and from the existence of a standard distinguished triangle associated to a closed inclusion.
\er

Let $i: Y\hookrightarrow Z$ be a closed inclusion of complex algebraic varieties
 with open complement $j: U=Z\backslash Y\hookrightarrow Z$. Then one has by \cite{Sa1}[(4.4.1)]   the following functorial distinguished triangle in $D^b\MHM(Z)$: 
\begin{equation}\label{tri}
 \begin{CD}
  j_!j^* @> ad_j >> id @> ad_i >> i_*i^* @> [1]>> \:,
 \end{CD}
\end{equation}
where the maps $ad$ are the adjunction maps, and $i_*=i_!$ since $i$ is proper.
In particular, we get the following {\em additivity relation} at the level of Grothendieck groups:
\be [\bQ^H_Z] =  [j_!\bQ^H_U] +  [i_!\bQ^H_Y] \in K_0(D^b\MHM(Z))=K_0(\MHM(Z)).
\ee
As a consequence, if $\cS=\{S\}$ is a complex algebraic stratification of $Z$ such that for any $S \in \cS$, $S$ is smooth and $\bar{S} \setminus S$ is a union of strata, then with $j_S:S \hookrightarrow Z$ denoting the corresponding inclusion map, we get:
\be [\bQ^H_Z] =\sum_{S \in \cS} [(j_S)_!\bQ^H_S]. 
\ee
If $f: Z \to X$ is a complex algebraic morphism, then we can apply $f_!$ to (\ref{tri}) to get another 
distinguished triangle
\begin{equation}\label{triangle2}
 \begin{CD}
  f_!j_!j^*\bQ^H_Z @> ad_j >> f_!\bQ^H_Z @> ad_i >> f_!i_!i^*\bQ^H_Z @> [1]>> \:.
 \end{CD}
\end{equation}
So, on the level of Grothendieck groups, we get the following additivity relation:
\begin{equation}\label{tri3}
[f_!\bQ^H_Z] = [(f\circ j)_!\bQ^H_U] + [(f\circ i)_!\bQ^H_Y] \in K_0(D^b\MHM(X))=K_0(\MHM(X))\:.
\end{equation}

Let $K_0(\var/X)$ be the motivic {\em relative Grothendieck group} of complex algebraic varieties over $X$,
i.e., the free abelian group generated by isomorphism classes $[f]:=[f: Z\to X]$ of morphisms $f$ to $X$,
divided out be the {\em additivity relation} 
$$[f]=[f\circ i] + [f\circ j]$$ 
for a closed inclusion
$i: Y\hookrightarrow Z$ with open complement $j: U=Z\backslash Y\hookrightarrow Z$. The  pushdown $f_!$, external product $\boxtimes$ and pullback $g^*$
for these relative Grothendieck groups 
are defined by composition, exterior product and resp. pullback of arrows. Then we get by (\ref{tri3}) the following result:
 
\bc\label{chiHdg}
There is a natural group homomorphism
$$\chi_{Hdg}: K_0(\var/X) \to K_0(\MHM(X)), \ \ \ [f:Z\to X]\mapsto [f_!\bQ^H_Z]\:,$$
which commutes with pushdown $f_!$, exterior product $\boxtimes$ and pullback $g^*$. 
\ec
The fact that $\chi_{Hdg}$ commutes
with external product $\boxtimes$ (or pullback $g^*$) follows from the corresponding K\"{u}nneth (or base change) theorem
for the functor 
$$f_!: D^b\MHM(Z)\to D^b\MHM(X)$$ 
(see \cite{Sa1}[(4.4.3)]).


\section{Hodge-Chern and Hirzebruch Classes of Singular Varieties}\label{HCH}

\subsection{Construction. Properties}

The construction of $K$-theoretical and resp. homology characteristic classes of mixed Hodge modules is based on the following result of Saito
(see \cite{Sa0}[Sec.2.3] and \cite{Sa3}[Sec.1] for the first part, and \cite{Sa1}[Sec.3.10, Prop.3.11]) for part two):

\bt[M. Saito]\label{grDR}
Let $X$ be a complex algebraic variety. Then there is a functor of triangulated categories
\begin{equation} Gr^F_p\DR: D^b\MHM(X) \to D^b_{\rm coh}(X)\end{equation}
commuting with proper pushforward,
and with $Gr^F_p\DR(\cM)=0$ for almost all $p$ and $\cM$ fixed,
where $D^b_{\rm coh}(X)$ is the bounded
derived category of sheaves of algebraic $\cO_X$-modules with coherent
cohomology sheaves. If $X$ is a (pure $n$-dimensional) complex algebraic manifold, then one has in addition the following:
\begin{enumerate}
\item Let $\cM\in \MHM(X)$ be a single mixed Hodge module. Then
$Gr^F_p\DR(\cM)$ is the corresponding complex associated to the de Rham complex of the underlying algebraic left $\cD$-module $\cM$ with its integrable connection $\nabla$:
$$ \begin{CD} \DR(\cM)=[\cM @> \nabla >> \cdots @> \nabla >> \cM\otimes_{\cO_X} \Omega^n_X]
\end{CD}$$ 
with $\cM$ in degree $-n$, filtered  by
$$ \begin{CD}
F_p \DR(\cM) =[F_p\cM @> \nabla >> \cdots  @> \nabla >> F_{p+n}\cM\otimes_{\cO_X} \Omega^n_X] \:.
\end{CD} $$
\item Let $\bar{X}$ be a smooth partial compactification of the complex algebraic manifold $X$ with complement $D$ a normal crossing divisor with smooth irreducible components, and with $j: X\to \bar{X}$ the open inclusion. 
Let $\bV=(L,F,W)$ be a good variation of mixed Hodge structures on $X$. 
Then
 the filtered de Rham complex 
$(\DR(j_*\bV),F)$ of  the shifted mixed Hodge module $j_*\bV\in \MHM(\bar{X})[-n]\subset D^b\MHM(\bar{X})$
is filtered quasi-isomorphic to the logarithmic de Rham complex 
$$\begin{CD} \DR_{\log}(\cL):=[\overline{\cL} @> \overline{\nabla} >> \cdots @> \overline{\nabla} >> \overline{\cL}\otimes_{\cO_{\bar{X}}} \Omega^n_{\bar{X}}(\log(D))] \end{CD}$$ with the increasing filtration $F_{-p}:=F^{p}$ ($p\in \bZ$) associated to the decreasing $F$-filtration
$$\begin{CD}
F^{p}\DR_{log}\left({\cL}\right) =
[F^{p}\overline{\cL} @> \overline{\nabla} >> \cdots @> \overline{\nabla} >> F^{p-n}\overline{\cL}\otimes_{\cO_{\bar{X}}} \Omega^n_{\bar{X}}(\log(D))] \:,
\end{CD}$$
where $\overline{\cL}$ is the canonical Deligne extension of $\cL:=L \otimes_{\bQ_X} \cO_X$.
 In particular, $Gr^F_{-p}\DR(j_*\bV)$ ($p\in \bZ$) is quasi-isomorphic to
$$\begin{CD}
 Gr_F^{p}\DR_{log}\left({\cL}\right) =
[Gr_F^{p}\overline{\cL} @> Gr\;\overline{\nabla} >> \cdots @> Gr\;\overline{\nabla} >> Gr_F^{p-n}\overline{\cL}\otimes_{\cO_{\bar{X}}} \Omega^n_{\bar{X}}(\log(D))] \:.
\end{CD}$$
 Similar considerations apply  to the filtered de Rham complex 
$(\DR(j_!\bV),F)$ of  the shifted mixed Hodge module $j_!\bV\in \MHM(\bar{X})[-n]\subset D^b\MHM(\bar{X})$, by considering instead the logarithmic de Rham complex  associated to the Deligne extension ${\overline \cL} \otimes \cO(-D)$ of $\cL$.
\end{enumerate}
\et

Note that the maps $Gr{\nabla}$ and $Gr\;\overline{\nabla}$ in the complexes $Gr_F^p\DR(\cM)$ and respectively $Gr_F^{p}\DR_{log}\left({\cL}\right)$ are $\cO$-linear.

The transformations
$Gr^F_p\DR$ ($p\in \bZ$) induce functors on the level of Grothendieck groups.
Therefore, if $K_0(X) \simeq K_0(D^b_{\rm coh}(X))$ denotes the
Grothendieck group of coherent {\em algebraic} $\cO_X$-sheaves on $X$, we get group
homomorphisms 
$$Gr^F_p\DR: K_0(\MHM(X))=K_0(D^b\MHM(X))\to K_0(D^b_{\rm coh}(X))\simeq K_0(X)\:.$$

\noindent{\bf Notation.} In the following we will use the notation $Gr^p_F$ in place of $Gr^F_{-p}$ corresponding to  the identification $F^{p}=F_{-p}$, which makes the transition between
the {\it increasing} filtration  $F_{-p}$ appearing in the $\cD$-module language and 
 the classical situation of the {\it decreasing} filtration $F^{p}$ coming from a good variation of mixed Hodge structures.

\bd\label{def12} For a complex algebraic variety $X$, the {\it Hodge-Chern class transformation} 
$$\DR_y:K_0(\MHM(X)) \to K_0(X) \otimes \bZ[y^{\pm 1}]$$
is defined by 
\be [\cM] \mapsto \DR_y([\cM]):=\sum_{i,p}\,(-1)^i\,\big[\cH^i\Gr_F^p\DR(\cM)\big]
\cdot (-y)^p \ee
(where the sign of the variable $y$ is chosen to fit with Hirzebruch's convention in (\ref{37}) below).
The  {\it un-normalized homology Hirzebruch class transformation} is defined as the composition \be\label{tr} T_{y*}:=td_* \circ \DR_y :K_0(\MHM(X)) \to H_*(X) \otimes \bQ[y,y^{-1}]\ee
where $td_*:K_0(X) \to H_*(X)\otimes \bQ$ is the Baum-Fulton-MacPherson Todd class transformation \cite{BFM}, which is linearly extended over $\bZ\big[y, y^{-1} \big]$. Here, we denote by $H_*(X)$ the even degree Borel-Moore homology $H_{2*}^{BM}(X)$ of $X$. 
The  {\it normalized homology Hirzebruch class transformation} is defined as the composition \be\label{trn}  \widehat{T}_{y*}:=td_{(1+y)*} \circ \DR_y :K_0(\MHM(X)) \to H_*(X)\otimes \bQ[y^{\pm 1},\frac{1}{y+1}]\ee
where $$td_{(1+y)*}:K_0(X)[y,y^{-1}]\to
H_{*}(X)\otimes \bQ\big[y,\h{$\frac{1}{y(y+1)}$}\big]$$
is the scalar extension of the Todd class transformation
$td_*$ together with the multiplication by
$(1+y)^{-k}$ on the degree $k$ component.
\ed

\br
By precomposing with the transformation $\chi_{Hdg}$ of Corollary \ref{chiHdg}, we get similar {\it motivic} transformations, Chern and resp. Hirzebruch, defined on the relative Grothendieck group of complex algebraic varieties. It is the (normalized) motivic Hirzebruch class transformation which unifies (in the sense of \cite{BSY}) the previously known characteristic class theories for singular varieties, namely, the (rational) Chern class transformation of MacPherson \cite{MP} for $y=-1$, the Todd class transformation of Baum-Fulton-MacPherson \cite{BFM} for $y=0$, and the $L$-class transformation of Cappell-Shaneson \cite{CS} for $y=1$,  thus answering positively an old question of MacPherson about the existence of such a unifying theory (cf. \cite{MP77,Y98}). For the rest of the paper, we mostly neglect the more topological $L$-class transformation  related to self-dual constructible sheaf complexes (e.g., intersection cohomology complexes), because the above-mentioned relation to Hirzebruch classes can only be defined at the motivic level, but not on the Grothendieck group of mixed Hodge modules.
\er

\br By \cite{Sc}[Prop.5.21], we have that
\be \widehat{T}_{y*}([\cM])\in H_*(X)\otimes\bQ[y,y^{-1}],\ee
and, moreover, the specialization at $y=-1$,\be \widehat{T}_{-1*}([\cM])=c_*([\rat(\cM)]) \in H_*(X) \otimes \bQ \ee 
equals the MacPherson-Chern class of the underlying constructible sheaf complex $\rat(\cM)$ (i.e., the MacPherson-Chern class of the constructible function defined by taking stalkwise the Euler characteristic).
\er

The {\it homology Hirzebruch characteristic classes} of a complex algebraic variety $X$,  denoted by $T_{y*}(X)$ and resp. $\widehat{T}_{y*}(X)$, are obtained by evaluating the Hirzebruch transformations of Definition \ref{def12} on the constant Hodge module $\cM=\bQ^H_{X}$. Moreover, we have that 
\be\label{101} T_{y*}(X), \ \widehat{T}_{y*}(X) \in H_*(X) \otimes \bQ[y].\ee

If  $X$ is smooth, then
\be\label{37}
\DR_y([\bQ^H_X])=\La_y[T^*_X],
\ee
where for a vector bundle $V$ on $X$ we set: $$\La_y[V]:=\sum_p[\La^pV]\,y^p.$$
Indeed, we have that
$$\DR(\bQ_X^H)=\DR(\cO_X)[-n]=\Omega_X^{\bullet},\q\h{with}\,\,\,
n:=\dim X,
$$
and the Hodge filtration $F^p$ on $\Omega_X^{\bullet}$ is given by the
truncation $\sig_{\ge p}$.
In particular, the following {\it normalization} property holds for the homology Hirzebruch classes: if $X$ is smooth, then
\be T_{y*}(X)=T_y^*(T_X) \cap [X]\:, \ \ \widehat{T}_{y*}(X)=\widehat{T}_y^*(T_X) \cap [X]\:,\ee
with $T_y^*(T_X)$ and $\widehat{T}_y^*(T_X)$ the un-normalized and resp. normalized versions of the 
cohomology Hirzebruch class of $X$ appearing in the generalized Riemann-Roch theorem, see \cite{H}. 
More precisely, we have un-normalized and resp. normalized power series
\be\label{ps} Q_y(\alpha):=\frac{\alpha(1+ye^{-\alpha})}{1-e^{-\alpha}} \ ,\,\,\,
\widehat{Q}_y(\alpha):=\frac{\alpha(1+ye^{-\alpha(1+y)})}{1-e^{-\alpha(1+y)}}
\in
\bQ[y][[\alpha]]\ee
with initial terms 
\be Q_y(0)=1+y, \q \widehat{Q}_y(0)=1\ee
so that for $X$ smooth we have 
\be T^*_y(T_X)=\prod_{i=1}^{\dim X}Q_y(\alpha_i)\in
H^{*}(X) \otimes\bQ[y] \ee
(and similarly for $\widehat{T}_y^*(T_X)$), with $\{\alpha_i\}$  the Chern roots of $T_X$. 
These two power series are related by the following relation
\be \widehat{Q}_y(\alpha)=(1+y)^{-1}\cdot Q_y(\alpha(1+y)),\ee
which explains the use of the normalized Todd class transformation $td_{(1+y)*}$ in the definition of $\widehat{T}_{y*}(X)$.

Note that for the values $y=-1$, $0$, $1$ of the parameter,  the class $\widehat{T}_y^*$ specializes to the total Chern class $c^*$, Todd class $td^*$, and $L$-polynomial $L^*$, respectively. Indeed, the power series $\widehat{Q}_y(\alpha)\in\bQ[y][[\alpha]]$ becomes respectively
$$1+\alpha,\q\alpha/(1 -e^{-\alpha}),\q\alpha/\tanh\alpha.$$
Moreover, by (\ref{101}) we are also allowed to specialize the parameter $y$ in the homology classes $ T_{y*}(X)$ and $\widehat{T}_{y*}(X)$
to the distinguished values $y=-1,0,1$. For example, for $y=-1$, we have the identification (cf. \cite{BSY})
\be\label{-1}\widehat{T}_{-1*}(X)=c_*(X) \otimes \bQ\ee
with the total (rational) Chern
class $c_*(X)$ of MacPherson \cite{MP}. 
Also, for a variety $X$ with at most {\it Du Bois singularities}
(e.g., rational singularities, such as toric varieties), we have by \cite{BSY} that
\be\label{y0}T_{0*}(X)=\widehat{T}_{0*}(X)=td_*([\mathcal{O}_X])=:td_*(X) \:,\ee
for  $td_*$ the
Baum-Fulton-MacPherson Todd class transformation \cite{BFM}. Indeed, in the language of mixed Hodge modules, $X$ has at most Du Bois singularities if, and only if, the canonical map 
$$\begin{CD} \cO_X @>\sim>> Gr^0_F\DR (\bQ^H_X) \in D^b_{\rm coh}(X) \end{CD}$$
is a quasi-isomorphism (see \cite{Sa3}[Cor.0.3]).
It is still only conjectured that if $X$ is a compact  complex algebraic variety which is also a rational homology manifold, then $\widehat{T}_{1*}(X)$ is the Milnor-Thom homology  $L$-class of $X$ (see \cite{BSY}). This conjecture is proved in \cite{CMSS2}[Cor.1.2] for projective varieties $X$ of the form $Y/G$, with $Y$ a projective manifold and $G$ a finite group of algebraic automorphisms of $Y$, see Remark \ref{Za}.\\

By the Riemann-Roch theorem of \cite{BFM}, the Todd class transformation  
$td_*$ commutes with the pushforward under proper morphisms, 
so the same is true for  the Hirzebruch transformations $T_{y*}$ and $\widehat{T}_{y*}$ of Definition \ref{def12}.
If we apply this observation to the constant map $k:X\to pt$ with $X$ compact, then the
pushforward for $H_*$ is identified with the degree map, and we have
$$K_0(pt)=\bZ,\q H_*(pt)=\bQ,\q\MHM(pt)=\MHS^p,$$
where $\MHS^p$ is, as before,  Deligne's category of graded-polarizable mixed $\bQ$-Hodge
structures.
By definition, we have for $H^{\bullet}\in D^b\MHS^p$ that
\be\label{316} T_{y*}(H^{\bullet})=\widehat{T}_{y*}(H^{\bullet})=\chi_y(H^{\bullet}):=
\sum_{j,p}\,(-1)^j\dim_{\bC}\Gr_F^pH_{\bC}^j\,(-y)^p.\ee
Hence, for $X$ compact and connected, the degree-zero part of $T_{y*}(X)$ or $\widehat{T}_{y*}(X)$ is identified with the {\it Hodge polynomial} 
\be \chi_y(X):=\chi_y\big(H^{\bullet}(X)\big),\ee
i.e.,
\be
 \chi_y(X)=\int_X  T_{y*}(X)=\int_X \widehat{T}_{y*}(X).
\ee

\br The Hodge polynomial of (\ref{316}) is just a special case of the {\it Hodge-Deligne polynomial} defined by taking into consideration both the Hodge and the weight filtration. However, simple examples show that there  is no characteristic class theory incorporating the weight filtration of mixed Hodge modules (e.g., see \cite{BSY}[Ex.5.1] and \cite{Sc}[p.428]).
\er

The Hodge-Chern and resp. Hirzebruch characteristic classes have good functorial properties, but are very difficult to compute in general. In the following, we use Theorem \ref{grDR}  to compute these characteristic classes in some simple examples. 

\bex\label{AS}
Let $\bar{X}$ be a smooth partial compactification of the complex algebraic manifold $X$ with complement $D$ a normal crossing divisor with smooth irreducible components, with $j: X\hookrightarrow \bar{X}$ the open inclusion. 
Let $\bV=(L,F,W)$ be a good variation of mixed Hodge structures on $X$. 
Then the filtered de Rham complex 
$(\DR(j_*\bV),F)$ of $j_*\bV\in \MHM(\bar{X})[-n]\subset D^b\MHM(\bar{X})$
is by Theorem \ref{grDR}(2) filtered quasi-isomorphic to the logarithmic de Rham complex $\DR_{log}(\cL)$ with the increasing filtration $F_{-p}:=F^{p}$ ($p\in \bZ$) associated to the decreasing $F$-filtration induced by Griffiths' transversality. Then
\begin{equation}\begin{split}\label{MHCj}
\DR_y([j_*\bV]) &= \sum_{i,p}\; (-1)^{i} [\cH^i ( Gr_F^{p} \DR_{log}(\cL) )] \cdot (-y)^p\\
&=\sum_{p}\; [Gr_F^{p} \DR_{log}(\cL) ] \cdot (-y)^p\\
&= \sum_{i,p}\; (-1)^{i} [ Gr_F^{p-i}(\overline{\cL}) \otimes_{\cO_{\bar{X}}} \Omega_{\bar{X}}^i(\log(D))] \cdot (-y)^p,\\
\end{split}\end{equation}
where the last equality uses the fact that the complex $Gr_F^{p} \DR_{log}(\cL) $ has coherent (locally free) objects, with $\cO_{\bar X}$-linear maps between them.
Let us now define $$Gr^y(Rj_*L):=\sum_p [Gr_F^{p}(\overline{\cL})] \cdot (-y)^p \in K^0({\bar X})[y^{\pm 1}],$$
with $K^0({\bar X})$ denoting the Grothendieck group of algebraic vector bundles on $\bar X$.
Therefore, (\ref{MHCj}) is equivalent to the formula:
\be \DR_y([j_*\bV])=Gr^y(Rj_*L) \cap  \left(\Lambda_y\left( \Omega_{\bar{X}}^1(\log(D))\right)\cap [\cO_{\bar{X}}] \right).
\ee
Similarly, by setting
$$Gr^y(j_!L):=\sum_p [\cO_{\bar X}(-D) \otimes Gr_F^{p}(\overline{\cL})] \cdot (-y)^p \in K^0({\bar X})[y^{\pm 1}],$$
we obtain the identity:
\be \DR_y([j_!\bV])=Gr^y(j_!L) \cap  \left(\Lambda_y\left( \Omega_{\bar{X}}^1(\log(D))\right)\cap [\cO_{\bar{X}}] \right).
\ee
Here, the pairing $$\cap:=\otimes_{\cO_{\bar X}}:K^0({\bar X}) \times K_0({\bar X}) \to K_0({\bar X})$$ is defined by taking the tensor product.
In particular, for $j=id: X\to X$ we get the following {\em Atiyah-Meyer type formula}
(compare \cite{CLMS, MSc}):
\begin{equation}\label{MHCy-fomula}
\DR_y([\bV]) = Gr^y(L) \cap \left(\Lambda_y(T^*_X)\cap [\cO_X] \right)\:.
\end{equation}
\eex

\bigskip

Let us now discuss the following {\it multiplicativity property} of the Hodge-Chern and resp. homology Hirzebruch class transformations:
\bp The Hodge-Chern class transformation
$\DR_y$ commutes with external products, i.e.:
\begin{equation}
\DR_y([\cM\boxtimes \cM'])=
\DR_y([\cM]\boxtimes [\cM'])= \DR_y([\cM]) \boxtimes \DR_y([\cM'])
\end{equation}
for $\cM\in D^b\MHM(Z)$ and $\cM'\in D^b\MHM(Z')$. A similar property holds for the (un-) normalized Hirzebruch class tranformations.
\ep
\noindent First note that since the Todd class transformation $td_*$ commutes with external products, it suffices to justify the first part of the above proposition (refering to the Hodge-Chern class).
For $X$ and $X'$ and their partial compactifications as in Example \ref{AS}, we have that:
$$\Omega_{\bar{X}\times \bar{X}'}^1(\log(D\times X'\cup X\times D')) \simeq \left(\Omega_{\bar{X}}^1(\log(D))\right) \boxtimes
\left( \Omega_{\bar{X'}}^1(\log(D'))\right).$$
Therefore, 
\begin{multline*} \Lambda_y\left( \Omega_{\bar{X}\times \bar{X}'}^1(\log(D\times X'\cup X\times D'))\right)\cap [\cO_{\bar{X}\times \bar{X}'}] = \\
\left(\Lambda_y\left( \Omega_{\bar{X}}^1(\log(D))\right)\cap [\cO_{\bar{X}}] \right)\boxtimes
\left(\Lambda_y\left( \Omega_{\bar{X'}}^1(\log(D'))\right)\cap [\cO_{\bar{X'}}] \right).\end{multline*}
Recall that the Grothendieck group $K_0(\MHM(Z))$ of mixed Hodge modules on the complex variety $Z$ is generated by classes of the form $f_*(j_*[\bV])$, with $f: \bar{X}\to Z$ proper and $X,\bar{X},\bV$ as
before. Finally one also has the multiplicativity
$$(f\times f')_*= f_*\boxtimes f'_*$$
for the pushforward for proper maps $f: \bar{X}\to Z$ and  $f': \bar{X}'\to Z'$
on the level of Grothendieck groups $K_0(\MHM(-))$, as well as for $K_0(-)\otimes\bZ[y^{\pm 1}]$. So the claim follows.\\

We conclude this section with a discussion of the following {\it additivity property} of 
the Hodge-Chern classes and resp. homology Hirzebruch classes (see \cite{MSS2}[Prop.5.1.2]):
\bp\label{p1}
Let $X$ be a complex algebraic variety, and fix $\cM \in D^b\MHM(X)$ with underlying $\bQ$-complex $K$.
Let $\cS=\{S\}$ be a complex algebraic stratification of $X$ such
that for any $S\in\cS$, $S$ is smooth, $\So\setminus S$ is a union of
strata, and the sheaves $\cH^iK|_S$ are local systems on $S$ for any $i$.
Let $j_S:S\hookrightarrow X$ denote the inclusion map. Then
\be 
\begin{split} [\cM] &=\sum_{S} [(j_S)_!(j_S)^*\cM] = \sum_{S,i}\,(-1)^i\,[(j_S)_!H^{i}
(j_S)^*\cM ] 
\\ &=\sum_{S,i}\,(-1)^i\,\big[(j_S)_!H^{i+\dim(S)} (j_S)^*\cM [-\dim(S)]\big] \in K_0(\MHM(X)), 
\end{split}
\ee
where $H^{i+\dim(S)}(j_S)^*\cM$ is a smooth mixed Hodge module on the stratum $S$ so that
$\cH^iK|_S \simeq {^p}\cH^{i+\dim(S)}(K|_S)[-\dim(S)]$ underlies a good variation of mixed Hodge structures.

 In particular, 
\be\label{0b}
\DR_y([\cM])=\sum_{S,i}\,(-1)^i\,\DR_y\big[(j_S)_!H^{i+\dim(S)} (j_S)^*\cM [-\dim(S)]\big],
\ee
and
\be\label{0}  T_{y*}(\cM)  =\sum_{S,i}\,(-1)^i\,T_{y*}\big((j_S)_!H^{i+\dim(S)} (j_S)^*\cM [-\dim(S)]\big).\ee
 A similar formula holds for the normalized Hirzebruch classes \tyh.
\ep
Moreover, the summands on the right-hand side of formulae (\ref{0b}) and resp. (\ref{0}) can be computed as in \cite{MSS2}[Sect.5.2] as follows. For any sets $A \subset B$, denote the inclusion map by $i_{A,B}$.
Let $\bV$ be a good variation of mixed Hodge structure on a
stratum $S$. 
Take a smooth partial compactification $i_{S,Z}:S\hookrightarrow Z$ such that
$D:=Z\setminus S$ is a divisor with simple normal crossings and,
moreover, $i_{S,\So}=\pi_Z\circ i_{S,Z}$ for a proper morphism
$\pi_Z:Z\to\So$.
Then in the notations of Example \ref{AS}, we get by functoriality the following result:
\bp\label{39} In the above notations, we have:
\be\label{00b} \DR_y(\big[ (i_{S,\So})_!\bV\big])=
(\pi_Z)_*\big[Gr^y((i_{S,Z})_!\bV) \cap  \Lambda_y\left( \Omega_{Z}^1(\log(D))\right].
\ee
In particular, if ${\overline \cL}$ denotes the canonical Deligne extension on $Z$ associated to  $\bV$, we obtain:
\be\label{00} T_{y*}\big((i_{S,\So})_!\bV\big)=
\sum_{p,q}(-1)^q(\pi_Z)_*td_{*}\big[\cO_Z(-D) \otimes \Gr_F^p{\overline \cL}
\otimes\Omega_Z^q(\log D)\big](-y)^{p+q}.
\ee
\ep

Let us now consider the situation of Proposition \ref{p1} in the special case when $\cM=\bQ^H_X$ is the constant Hodge module. Then the following additivity holds:
\be  [\bQ_X^H]=\sum_{S}\,\big[(j_S)_!\bQ^H_S\big]=\sum_{S}\,(i_{\So,X})_*\big[(i_{S,\So})_!\bQ^H_S\big] \in K_0(\MHM(X)),
\ee
where we use the factorization $j_S=i_{\So,X} \circ i_{S,\So}$.
Hence the formulae (\ref{0b}) and (\ref{0}) yield by functoriality the following:
\be\label{000b} \DR_y([\bQ_X^H])=
\sum_{S}\,(i_{\So,X})_*\DR_y\big[(i_{S,\So})_!\bQ^H_S\big], \ee
and
\be\label{000} T_{y*}(X)
=\sum_{S}\,(i_{\So,X})_*T_{y*}\big((i_{S,\So})_!\bQ^H_S\big). \ee
Moreover, in the notations of Proposition \ref{39}, we get by (\ref{00b}):
\be\label{0000b}
 \DR_y(\big[(i_{S,\So})_!\bQ^H_S\big]) =(\pi_Z)_*\big[\cO_{Z}(-D)
\otimes \Lambda_y \Omega_Z^1(\log D)\big] \in K_0({\bar S})[y],
\ee
and
\be\label{0000}
\begin{split} T_{y*}\big((i_{S,\So})_!\bQ^H_S\big) &=
\sum_{q\geq 0} (\pi_Z)_*td_{*}\big[\cO_{Z}(-D)
\otimes\Omega_Z^q(\log D)\big]y^q \\ &=(\pi_Z)_*td_{*}\big[\cO_{Z}(-D)
\otimes \Lambda_y \Omega_Z^1(\log D)\big] \in H_*({\bar S}) \otimes \bQ[y],
\end{split}
\ee
where for a vector bundle $V$ we set as before: $\Lambda_y[V]:=\sum_p[\Lambda^pV]\,y^p$.


\subsection{Application: Characteristic Classes of Toric Varieties}\label{Toric}

Let $X_{\Sig}$ be a toric variety of dimension $n$ corresponding to the fan $\Sigma$, and with torus $T:=(\bC^*)^{n}$; see \cite{CLS} for more details on toric varieties. Then $X_{\Sig}$ is stratified by the orbits of the torus action. More precisely, by the orbit-cone correspondence, to a $k$-dimensional cone $\sigma \in \Sig$ there corresponds an $(n-k)$-dimensional torus-orbit \os $ \ \cong (\bC^*)^{n-k}$.  The closure \vs \ of \os \ (which is the same in both classical and Zariski topology) is itself a toric variety, and a $T$-invariant subvariety of $X_{\Sig}$. In particular, we can apply the results of the previous section in the setting of toric varieties.

By using toric geometry, formula (\ref{0000b}) adapted to the notations of this section yields the following 
(see \cite{MS13}[Prop.3.2]):
\bp\label{p3} For any cone $\sig \in \Sig$, with corresponding orbit \os \ and inclusion $i_{\sig} : O_{\sig} \hookrightarrow \overline{O}_{\sig}=V_{\sig}$, we have:
\be
\DR_y(\big[(i_{\sig})_!\bQ^H_{O_{\sig}}\big])=(1+y)^{\dim(O_{\sig})} \cdot [\omega_{V_{\sig}}],
\ee
where $\omega_{V_{\sig}}$ is the canonical sheaf on the toric variety $V_{\sig}$.
\ep

The main result of this section is now a direct consequence of Propositions \ref{p3} and additivity:
\bt\label{t1}
Let $X_{\Sig}$ be the toric variety defined by the fan $\Sigma$. For each cone $\sig \in \Sig$ with corresponding torus-orbit $O_{\sig}$, denote by $k_{\sig} : V_{\sig} \hookrightarrow X_{\Sig}$ the inclusion of the orbit closure.   Then the Hodge-Chern class of $X_{\Sig}$ is computed by the formula:
\be
\DR_y(X_{\Sig})=\sum_{\sig \in \Sig} \, (1+y)^{\dim(O_{\sig})} \cdot (k_{\sig})_*[\omega_{V_{\sig}}].
\ee
Similarly, the 
un-normalized homology Hirzebruch class $T_{y*}(X_{\Sig})$ is computed by:
\be\label{7} 
T_{y*}(X_{\Sig})=\sum_{\sig \in \Sig}\, (1+y)^{\dim(O_{\sig})} \cdot (k_{\sig})_*td_*([\omega_{V_{\sig}}]).
\ee
And after normalizing, the corresponding  homology Hirzebruch class $\widehat{T}_{y*}(X_{\Sig})$ is computed by the following formula:
\be\label{8} 
\widehat{T}_{y*}(X_{\Sig})=\sum_{\sig, k}\, (1+y)^{\dim(O_{\sig}) -k} \cdot (k_{\sig})_*td_k([\omega_{V_{\sig}}]).
\ee
\et

Recall that by making $y=-1$ in the normalized Hirzebruch class $\widehat{T}_{y*}$ one gets the (rational) MacPherson Chern class $c_*$. Moreover, since toric varieties have only  rational (hence Du Bois) singularities, making $y=0$  in the Hirzebruch classes yields the Todd class $td_*$. 
So we get as a corollary of Theorem \ref{t1} the following result:
\bc 
The (rational) MacPherson-Chern class $c_*(X_{\Sig})$ of a toric variety $X_{\Sig}$ is computed by {\rm Ehler's formula}:
\be\label{E}
c_*(X_{\Sig})=\sum_{\sig \in \Sig}\,  (k_{\sig})_*td_{\dim(O_{\sig})}(V_{\sig})=\sum_{\sig \in \Sig}\,  (k_{\sig})_*([V_{\sig}]).
\ee
The Todd class $td_*(X_{\Sig})$ is computed by the additive formula:
\be\label{T}
td_*(X_{\Sig})=\sum_{\sig \in \Sig}\, (k_{\sig})_*td_*([\omega_{V_{\sig}}]).
\ee
\ec

\br The results of this section hold more generally for torus-invariant closed algebraic subsets of  $X_{\Sig}$ which are known to also have Du Bois singularities. For more applications and examples, e.g., generalized Pick-type formulae for full-dimensional lattice polytopes, see \cite{MS13}.
\er


\section{Hirzebruch-Milnor Classes}\label{HM}

In this section, we explain how to compute the homology Hirzebruch classes of globally defined hypersurfaces in an algebraic manifold (but see also \cite{MSS2} for the global complete intersection case). The key technical result needed for this calculation is the {\it specialization property} of the Hirzebruch class transformation, obtained by the second author in \cite{Sb}. We first explain this result in Section \ref{specv} after some relevant background is introduced in Section \ref{v}. We complete our calculation of Hirzebruch classes of hypersurfaces  in Section \ref{applic}.

\subsection{Motivation} 
Let $X \overset{i}{\hookrightarrow} M$ be the inclusion of an algebraic hypersurface $X$  in a  complex algebraic manifold $M$ (or more generally the inclusion of a local complete intersection). Then the normal cone $N_XM$ is a complex algebraic vector bundle $N_XM\to X$
over $X$, called the normal bundle of $X$ in $M$. The {\it virtual tangent bundle} of $X$, that is, 
\begin{equation}T^{{\rm vir}}_X:=[i^*T_M-N_XM] \in K^0(X),\end{equation} is independent of the embedding in $M$, so it is a well-defined element in the Grothendieck group $K^0(X)$ of algebraic vector bundles on $X$. Of course 
$$T^{{\rm vir}}_X=[T_X] \in K^0(X),$$
in case $X$ is a smooth algebraic submanifold.
Let $cl^*$ denote a multiplicative characteristic class theory of complex algebraic vector bundles, i.e., a natural transformation (with  $R$ a commutative ring with unit)
$$cl^*: \left(K^0(X),\oplus\right)\to \left(H^*(X)\otimes R,\cup\right)\:,$$
from the Grothendieck group $K^0(X)$ of complex algebraic vector bundles to a suitable cohomology theory $H^*(X)$  with a cup-product $\cup$, e.g., $H^{2*}(X;\bZ)$.
Then one can associate to $X$ an {\it intrinsic}  homology class (i.e., independent of the embedding $X\hookrightarrow M$) defined as follows:
\begin{equation}
cl_*^{\rm vir}(X):=cl^*(T^{{\rm vir}}_X) \cap [X] \in H_*(X)\otimes R\:.
\end{equation} 
Here $[X]\in H_*(X)$ is the {\it fundamental class} of $X$ in a suitable homology theory  $H_*(X)$ (e.g., the  Borel-Moore homology $H_{2*}^{BM}(X)$).\\

Assume, moreover, that there is a homology characteristic class theory $cl_*(-)$ for complex algebraic varieties, functorial for proper morphisms, obeying the normalization condition that for $X$ smooth  $cl_*(X)$ is the Poincar\'e dual of $cl^*(T_X)$ (justifying the notation $cl_*$).  If $X$ is smooth,  then clearly we have that  
$$cl_*^{\rm vir}(X)=cl^*(T_X) \cap [X]=cl_*(X) \:.$$ 
However, if $X$ is singular, the difference between  the homology classes $cl_*^{\rm vir}(X)$ and $cl_*(X)$ depends in general on the singularities of $X$. This motivates the following {problem}:\\ ``{\it Describe the difference class $cl_*^{\rm vir}(X)-cl_*(X)$ in terms of the geometry of the singular locus of $X$}.

This problem is usually studied in order to understand the complicated homology classes $cl_*(X)$ in terms of the simpler virtual classes $cl_*^{\rm vir}(X)$ and these difference terms measuring the complexity of singularities of $X$. The strata of the singular locus have a rich geometry, beginning with generalizations of knots which describe their local link pairs. This ``normal data", encoded in algebraic geometric terms via, e.g., the mixed Hodge structures on the (cohomology of the) corresponding Milnor fibers,  will play a fundamental role in our study of characteristic classes of hypersurfaces.\\

There are a few instances in the literature where, for the appropriate choice of $cl^*$ and $cl_*$, this problem has been solved. The first example was for the Todd classes $td^*$, and $td_{*}(X):=td_{*}([\cO_{X}])$, respectively, with
$$td_{*}: K_{0}(X) \to H_{*}(X)\otimes \bQ$$ 
the Todd class transformation in the singular Riemann-Roch theorem of
Baum-Fulton-MacPherson \cite{BFM}. Here $[\cO_X] \in K_0(X)$ the class of the structure sheaf.
By a result of Verdier \cite{V0, Fu}, $td_*$ commutes with the corresponding Gysin homomorphisms for the regular embedding $i: X\hookrightarrow M$.
This can be used to show that 
 $$td_*^{\rm vir}(X):=td^*(T^{{\rm vir}}_X) \cap [X]= td_*(X)$$ 
 equals the Baum-Fulton-MacPherson Todd class $td_*(X)$ of $X$ (\cite{V0, Fu}).\\

If $cl^*=c^*$ is the total Chern class in cohomology, the problem amounts to comparing the {\it Fulton-Johnson class} $c^{FJ}_*(X):=c_*^{vir}(X)$ (e.g., see \cite{Fu,FJ}) with the  homology Chern class $c_*(X)$ of MacPherson \cite{MP}. Here $c_*(X):=c_*(1_X)=c_*(\bQ_X)$, with
$$c_*: K_0(D^b_c(X)) \to F(X)\to H_*(X)$$
the functorial Chern class transformation of MacPherson \cite{MP}, defined on the group
$F(X)$ of complex algebraically constructible functions. To emphasize the analogy with the Grothendieck group of mixed Hodge modules (as used in the subsequent sections), we work here with the Grothendieck group of constructible (resp. perverse) sheaf complexes. 
The difference between the two classes $c_*^{vir}(X)$ and $c_*(X)$ is measured by the so-called {\it Milnor class}, $\cM_*(X)$. This is a homology class supported on the singular locus of $X$, and in the case of a global hypersurface $X=\{f=0\}$ it was computed in \cite{PP}  as a weighted sum in the Chern-MacPherson classes of closures of singular strata of $X$, the weights depending only on the normal information to the strata. For example, if $X$ has only isolated singularities, the Milnor class equals (up to a sign) the sum of the Milnor numbers attached to the singular points, which also explains the terminology:
\begin{equation}
\cM_*(X) =\sum_{x \in X_{\rm sing}} \;\chi\left(\tilde{H}^*(F_x;\bQ)\right)\:,
\end{equation}
where  $F_x$ is the local Milnor fiber of the isolated hypersurface singularity $(X,x)$.
More generally, Verdier \cite{V} proved the following  specialization result for the MacPherson-Chern class transformation (which holds even more generally for every constructible sheaf complex on $M$):
\be\label{Vc} c_*(\psi_f \bQ_M)=i^!c_*(\bQ_M),
\ee
where $i^!:H_*(M) \to H_{*-1}(X)$ is the homological Gysin map.
This result was used in  \cite{PP,Sch} for computing the   
(localized) Milnor class $\cM_*(X)$ of a global hypersurface $X=\{f=0\}$
in terms of the vanishing cycles of $f:M \to \bC$:
\begin{equation}\label{Milnor-van}
\cM_*(X)=c_*(\varphi_f(\bQ_M)) \in H_*(X_{\rm sing})\:,
\end{equation}
with the support of the (shifted) perverse complex $\varphi_f(\bQ_M)$ being contained in the singular locus $X_{\rm sing}$ of $X$.

\br For a more topological example concerning the Goresky-MacPherson $L$-class $L_*(X)$ (\cite{GM1}) for $X$ a {\it compact} complex hypersurface, see \cite{CS0,CS}.\er

A main goal here is to explain the (unifying) case when $cl^*=\widehat{T}_y^*$ is the (total) normalized cohomology Hirzebruch class of the generalized Hirzebruch-Riemann-Roch theorem \cite{H}.
The aim is to show that the results stated above are part of a more general philosophy, derived from comparing the intrinsic homology class (with polynomial coefficients) \begin{equation}{\widehat T}_{y*}^{\rm vir}(X):=\widehat{T}_y^*(T^{{\rm vir}}_X) \cap [X] \in H_*(X)\otimes \bQ[y]\end{equation}  
with the homology Hirzebruch class $\widehat{T}_{y*}(X)$ of \cite{BSY}.  This approach is motivated by  the fact  that, as already mentioned,  the $L$-class $L^*$, the Todd class $td^*$ and the Chern class $c^*$, respectively,  are all  suitable specializations (for $y=1,0,-1$, respectively) of the Hirzebruch class $\widehat{T}_y^*$; see \cite{H}.  In order to achieve our goal, we need to adapt Verdier's specialization result (\ref{Vc}) to the normalized Hirzebruch class transformation. For this, we first need to recall Saito's definition of nearby and vanishing cycles of mixed Hodge modules in terms of the $V$-filtration for the underlying filtered $\cD$-modules.


\subsection{$V$-filtration. Nearby and Vanishing Cycles}\label{v}
Let $f:M \to \bC$ be an algebraic function defined on a complex algebraic manifold $M$, with $X:=\{f=0\}$ a hypersurface in $M$. Consider the graph embedding $i':M \to M':=M \times \bC$, with $t=pr_2:M' \to \bC$ the projection onto the second factor. Note that $t$ is a smooth morphism, with $f=t \circ i'$.

Let ${\mathcal I}\subset {\mathcal O}_{M'}$ be the ideal sheaf defining the smooth hypersurface $\{t=0\} \simeq M$, i.e., the sheaf of functions vanishing along $M$. Then
the {\em increasing $V$-filtration of Malgrange-Kashiwara} on the algebraic coherent sheaf ${\mathcal D}_{M'}$ with respect to the smooth hypersurface $M\subset M'$ is defined for $k\in \bZ$ by
$$V_k{\mathcal D}_{M'}:=\{P\in {\mathcal D}_{M'}|\; P({\mathcal I}^{j+k})\subset {\mathcal I}^j
\quad \text{for all $j\in \bZ$}\}\:.$$
Here ${\mathcal I}^j:={\mathcal O_{M'}}$ for $j<0$. Note that
$$\bigcap_{k\in \bZ}\:V_k{\mathcal D}_{M'}=\{0\} \quad \text{and} \quad
\bigcup_{k\in \bZ}\:V_k{\mathcal D}_{M'}={\mathcal D}_{M'}\:.$$
Moreover,
$V_k{\mathcal D}_{M'}|_{\{t\neq 0\}}= {\mathcal D}_{M'}|_{\{t\neq 0\}}$ for all $k\in \bZ$,
so that $Gr^V_k{\mathcal D}_{M'}$ is supported on $M$. By definition, one has
$$t\in  V_{-1}{\mathcal D}_{M'}, \:\partial_t \in  V_{1}{\mathcal D}_{M'}
\quad \text{ and} \quad
\partial_t t = 1+t\partial_t \in  V_{0}{\mathcal D}_{M'} \:.$$ 
Also, for the sheaf 
${\mathcal D}_{M'/\bC}$
of {\em relative differential operators along the fibers of $t$} we have: 
$${\mathcal D}_{M'/\bC}\subset  V_{0}{\mathcal D}_{M'}
\quad \text{and} \quad Gr^V_0{\mathcal D}_{M'}|_M={\mathcal D}_{M}[\partial_t t] \:.$$

Let $\cM \in \MHM(M)$ be a given mixed Hodge module, where as before we use the same symbol for the underlying (filtered) holonomic (left) $\cD$-module. Then the pushforward $\cD_{M'}$-module $\cM':=i'_{*}\cM$ on $M'$  admits a unique increasing {\it canonical $V$-filtration} satisfying
 the following properties:
\begin{enumerate}
\item it a discrete, rationally indexed filtration.
\item $\bigcup_{\alpha} V_{\alpha} {\mathcal M'} = {\mathcal M'}$, and each $V_{\alpha} {\mathcal M'}$ is a coherent $V_{0}{\mathcal D}_{M'}$-module.
\item $(V_k{\mathcal D}_{M'}) (V_{\alpha} {\mathcal M'}) \subset V_{\alpha+k} {\mathcal M}$ for all $\alpha\in \bQ, k\in \bZ$, \\
and $t(V_{\alpha} {\mathcal M'}) = V_{\alpha-1} {\mathcal M'}$ for all $\alpha <0$.
\item $\partial_t t+\alpha$ is {\em nilpotent} on $Gr^V_{\alpha}{\mathcal M'}:= 
V_{\alpha} {\mathcal M'}/V_{<\alpha} {\mathcal M'}$, 
 with $V_{<\alpha} {\mathcal M'}:=\bigcup_{\beta<\alpha} V_{\beta} {\mathcal M'}$.
\end{enumerate}

The above properties are referred to as the {\it specializability} of $\cM'$ along $\{t=0\}$ (resp. of $\cM$ along $\{f=0\}$), where the rationality of the filtration corresponds to a quasi-unipotent  property.
These properties also imply that $t: Gr^V_{\alpha}{\mathcal M'}\to Gr^V_{\alpha-1}{\mathcal M'}$ is bijective for all
$\alpha \neq 0$, and $\partial_t: Gr^V_{\alpha}{\mathcal M'}\to Gr^V_{\alpha+1}{\mathcal M'}$ is bijective for all
$\alpha \neq -1$. 
Moreover,
\be\label{star} t\cdot: V_{\alpha} {\mathcal M'}\to V_{\alpha-1} {\mathcal M'}\quad
\text{is injective for all $\alpha<0$.}
\ee
Finally, all $Gr^V_{\alpha}{\mathcal M'}|_M$ are holonomic  left ${\mathcal D}_{M}$-modules.

\br
Here we work with an {\em increasing} $V$-filtration for {\em left} ${\mathcal D}$-modules,
so that $\partial_t t+\alpha$ is {\em nilpotent} on $Gr^V_{\alpha}{\mathcal M'}$.
\begin{enumerate}
\item This notion of $V$-filtration is compatible with the corresponding notion for the associated analytic $\cD$-modules. 
\item If we switch to the corresponding {\em right} ${\mathcal D}$-module 
$\omega_{M'}\otimes {\mathcal M'}$,
with the induced $V$-filtration $\omega_{M'}\otimes V_{\alpha}{\mathcal M'}$, then
$t\partial_t -\alpha$ is nilpotent on $Gr^V_{\alpha}(\omega_{M'}\otimes {\mathcal M'})$
(fitting with the convention of \cite{Sa0}[Def. 3.1.1]).
\item If we use the corresponding {\em decreasing} $V$-filtration $V^{\alpha}:=V_{-\alpha-1}$
for {\em left} ${\mathcal D}$-modules, then $t\partial_t -\alpha$ is nilpotent on $Gr_V^{\alpha}{\mathcal M'}$
(fitting with the convention of \cite{Sa0}[Introduction, p.851]).
\end{enumerate}
\er

By Saito's theory, we have the following compatibilities between the (Hodge) filtration $F$ and the $V$-filtration of the $\cD$-module underlying the mixed Hodge module $\cM'$:
\begin{enumerate}
\item[(s1)]  The induced $F$-filtration on $Gr^V_{\alpha} {\mathcal M'}$ is good for all $-1\leq \alpha \leq 0$.
\item[(s2)] $t(F_p V_{\alpha} {\mathcal M'}) = F_p V_{\alpha-1} {\mathcal M'}$ for all $\alpha <0$ and $p\in \bZ$.
\item[(s3)] $\partial_t (F_p Gr^V_{\alpha} {\mathcal M'}) = F_{p+1} Gr^V_{\alpha+1} {\mathcal M'}$ for all $\alpha >-1$ and $p\in \bZ$.
\end{enumerate}
The first property above is called the {\it regularity of $\cM$ with respect to $f$}, while the last two account for the {\it strict specializability of $\cM$ with respect to $f$}.
Here the {\em induced} $F$-filtration on $V_{\alpha}{\mathcal M'}$ resp. $Gr^V_{\alpha} {\mathcal M'}$ are given by
$$F_pV_{\alpha}{\mathcal M'}:= F_p{\mathcal M'}\cap  V_{\alpha}{\mathcal M'}$$ 
resp.
$$F_p Gr^V_{\alpha} {\mathcal M'}:=\im\left(F_pV_{\alpha}{\mathcal M'}\to Gr^V_{\alpha} {\mathcal M'}\right) \simeq F_pV_{\alpha}{\mathcal M'}/ F_pV_{<\alpha}{\mathcal M'} \:.$$

\br The above compatibility properties for algebraic mixed Hodge modules can be justified as follows:
\begin{enumerate}
\item By Saito's work \cite{Sa0}[Def.5.1.6], the underlying filtered $\cD$-module of an {\em analytic pure Hodge module} on $M'$ is {\em strictly specializable and  quasi-unipotent} along $M$.
\item In the complex algebraic context a  {\it pure} Hodge module of \cite{Sa0} is by definition ``extendable and 
quasi-unipotent'' at infinity, so that one can assume $M$ is compact, and then the properties
(s1-s3) in the analytic context of \cite{Sa0} imply by GAGA and flatness of
${\mathcal O}^{an}$ over ${\mathcal O}$ the same properties for the underlying algebraic filtrations.
\item The underlying filtered $\cD$-module of an {\em algebraic mixed Hodge module} on $M'$  is {\em strictly specializable and  quasi-unipotent} along $M$, since it is a finite successive extension of pure Hodge modules (by the weight filtration). The corresponding short exact extension sequences of mixed Hodge modules
are strict with respect to the (Hodge) filtration $F$ (cf. \cite{Sa0}[Lem.5]). 
Since the canonical $V$-filtration behaves well under extensions (cf. \cite{Sa1}[Cor.3.1.5]), one easily gets the
quasi-unipotence and the condition (s1-s3) above by induction (using also the properties of the underlying $V$-filtrations). Compare also with condition (2.2.1) in \cite{Sa1}[p.237].

\end{enumerate}
\er

The graded pieces $Gr^V_{\alpha}$ of the $V$-filtration are used in the $\cD$-module context for the definition of the exact nearby and resp. vanishing cycle functors  $$\psi_f^H, \varphi_f^H: \MHM(M) \to \MHM(X)$$ on the level of mixed Hodge modules.
If $${\rat}: \MHM(-) \to {\rm Perv}(\bQ_{-})$$ is the forgetful functor assigning to a mixed Hodge module the underlying $\bQ$-perverse sheaf, then $${\rat} \circ \psi^H_f={{^p}\psi_f} \circ {\rat}$$ and similarly for $\varphi_f^H$. Here ${{^p}\psi_f}:=\psi_f[-1]$ is a shift of Deligne's nearby cycle functor, and similarly for ${{^p}\varphi_f}$. 
So the shifted transformations $\psi'^H_f:=\psi^H_f[1]$ and $\varphi'^H_f:=\varphi^H_f[1]$
correspond under ${\rat}$ to the usual nearby and vanishing cycle functors. Recall that there exist decompositions $$\psi_f^H=\psi_{f,1}^H \oplus \psi_{f,\neq 1}^H$$ (and similarly for $\varphi_f^H$) into unipotent and non-unipotent part, with $\psi_{f,\neq 1}^H=\varphi_{f,\neq 1}^H$.

More precisely, for $\cM \in \MHM(M)$ (with $\cM':=i'_*\cM$), one has the following definition:
\bd \label{psi}
The underlying filtered $\cD$-module $\psi_f((\cM,F))$ of the {\it nearby cycle} mixed Hodge module $\psi_f^H(\cM)$ is defined by:
\begin{equation}
\psi_f\left(({\mathcal M},F)\right):= \bigoplus_{-1\leq \alpha <0} (Gr^V_{\alpha}{\mathcal M'}|_M,F) \:.
\end{equation}
The summand for $\alpha=-1$ (resp. $\alpha \neq -1$) corresponds to the (non-)unipotent nearby cycles $\psi_{f,1}^H(\cM)$ (resp. $\psi_{f,\neq 1}^H(\cM)$).
Similarly, the underlying filtered $\cD$-module of unipotent {\it vanishing cycles} mixed Hodge module $\varphi^H_{f,1}(\cM)$ is given by:
\be
\varphi_{f,1}\left((\cM,F)\right):=(Gr^V_{0}{\mathcal M'}|_M,F[-1]),
\ee
 with the shifted filtration defined as $(F[k])_i:=F_{i-k}$. \ed

\br Note that $\psi_f\left(({\mathcal M},F)\right)$ and $\varphi_{f,1}\left((\cM,F)\right)$ are filtered $\cD$-modules on the smooth hypersurface $M=\{t=0\}$ with support on $X$.

Our definition of the induced filtration fits with \cite{Sa0}[Introduction, p.851],
since we use {\em left} ${\mathcal D}$-modules. For the corresponding {\em right}
${\mathcal D}$-modules one has to shift these induced $F$-filtrations by $[+1]$
(see \cite{Sa0}[(5.1.3.3) on p.953]). This corrects the different switching by
$$\otimes (\omega_{M'},F) \quad \text{ or}  \quad \otimes (\omega_{M},F)$$ 
from filtered left ${\mathcal D}$-modules
to filtered right ${\mathcal D}$-modules on $M'$ or $M$, with $F$ the trivial filtration
such that $Gr^F_{-k}(-)= 0$ for $k\neq$ the dimension of the ambient manifold $M'$ or $M$.
\er


\subsection{Specialization of Hodge-Chern and Hirzebruch Classes}\label{specv}

In this section we discuss a Hodge-theoretic generalization of Verdier's result \cite{V} on the specialization of the MacPherson Chern class transformation.

Before formulating our main result, let us explain a motivating example in the case
when $i: X:=\{f=0\}\to M$ is a codimension one inclusion of complex algebraic {\it manifolds} (with $m:=\dim(M)$) and 
$$\bV^H={\mathcal M} \in \MHM(M)[-m]\subset D^b\MHM(M)$$ 
a (shifted) smooth mixed Hodge module corresponding to a good variation $\bV$ 
of mixed Hodge structures on $M$. Then by (\ref{MHCy-fomula}), we have an Atiyah-Meyer type formula
\begin{equation}
\DR_y([\bV]) = Gr^y(L) \cap \left(\Lambda_y(T^*_M)\cap [\cO_M] \right)\:,
\end{equation}
with $L$ denoting the underlying local system of the variation.
Then one gets
\begin{equation}\begin{split}
i^!\DR_y([{\bV}])&=i^*\left(Gr^y(L)\cup \Lambda_y(T^*_M)\right)
\cap i^!([{\mathcal O_M}])\\
&=\left( Gr^y(i^*L)\cup   \Lambda_y(i^*T^*_M) \right)\cap [{\mathcal O_{X}}]  \\
&= \Lambda_y(N^*_{X}M) \cap \DR_y([i^*{\bV}])\\
&= (1+y)\cdot \DR_y(i^*[{\bV}])\:.
\end{split}\end{equation}
Here we use the {\em multiplicativity} of $\Lambda_y(-)$ with respect to the short exact sequence of vector bundles
$$0\to N^*_{X}M \to T^*_M|_X \to T^*_X\to 0 \:,$$
and the triviality of the conormal bundle $N^*_{X}M$ (coming from the section $df$) so that
$$\Lambda_y(N^*_{X}M )= (1+y)\cdot [{\mathcal O_{X}}]\in K^0(X)[y]\:.$$
Since in this special case there are no vanishing cycles, $\varphi'^H_f({\bV})=0$, so 
$i^*\bV \simeq \psi'^H_f({\bV}).$
Therefore, we get:
\be\label{exspec} i^!\DR_y([{\bV}])=(1+y)\cdot \DR_y\left([\psi'^H_f({\bV})]\right)\:.\ee
As it will be explained below,  this formula (\ref{exspec}) holds in the general case of a (possibly singular) hypersurface $X=\{f=0\}$ of codimension one in the algebraic manifold $M$.  More precisely, by using the induced transformations 
$$\psi_f^H, \varphi_f^H: K_0(\MHM(M)) \to K_0(\MHM(X))$$
on the Grothendieck groups of mixed Hodge modules,  
we have the following counterpart of Verdier's
specialization result (see \cite{Sc}): 
\bt \label{motthm}
The Hodge-Chern class transformation $\DR_y$ commutes with specialization in the following sense:
\begin{equation}\label{spdr}
(1+y)\cdot \DR_y(\; \psi'^H_f (-)\;) = 
-(1+y)\cdot \DR_y(\; \psi^H_f (-)\;) = i^!\DR_y(-)   
\end{equation} 
as transformations $K_0(\MHM(M))\to K_0(X)[y,y^{-1}]$.
\et 
As an immediate corollary, by Verdier's specialization result for the Todd class \cite{V, Fu}, we have the following: 
\bc \label{main-sch}
The normalized Hirzebruch class transformation ${\widehat T}_{y*}$ commutes with specialization, that is: 
\begin{equation}\label{sp}
{\widehat T}_{y*}(\psi'^H_f(-))=i^! {\widehat T}_{y*}(-): \:K_0(\MHM(M)) \to
 H_*(X) \otimes \bQ[y,y^{-1}] \:.
\end{equation}
\ec
Note that the factor $(1+y)$ of equation (\ref{spdr}) disappears since $i^!$ sends $H_*(M)$ to $H_{*-1}(X)$.\\

We conclude this section with a sketch of proof of Theorem \ref{motthm}, see \cite{Sc} for complete details. Let us fix a mixed Hodge module $\cM$ on $M$. Then the proof uses only the underlying filtered $\cD$-module $(\cM,F)$ and its nearby cycles defined in terms of the $V$-filtration of $\cM':=i'_*\cM$ along the function $t=pr_2:M':=M \times \bC \to \bC$. Here $i':M \to M'$ is the graph embedding, with $f=t \circ i'$.  Note that, by construction, $$i^!\DR_y([\cM]) = i^!\DR_y([\cM']) ,$$ where on the right hand side $i$ denotes the hypersurface inclusion $(M=\{t=0\},X) \hookrightarrow (M',i'(M))$. So in the following we work only with the mixed Hodge module $\cM'$ on $M'$, with support in $i'(M) \simeq M$. 

On $M'$ we have the splitting: 
$$\Omega^1_{M'} = \Omega^1_{M'/\bC} \oplus  \Omega^1_{\bC}\:,$$
with $\Omega^1_{M'/\bC}$ the relative $1$-forms with respect to the submersion $t$
and $\Omega^1_{\bC} \simeq {\mathcal O}_{M'}dt$.
So we get:
$$\Omega^k_{M'} = \Omega^k_{M'/\bC} \oplus \Omega^1_{\bC}\otimes \Omega^{k-1}_{M'/\bC}\:.$$
This induces a splitting of the de Rham complex $\DR({\mathcal M'})$ as a double complex
$$\begin{CD} \DR_{/\bC}({\mathcal M'}) @=
[  \Omega^k_{M'/\bC} \otimes {\mathcal M'} @> \nabla_{/\bC} >>   \Omega^{k+1}_{M'/\bC} \otimes {\mathcal M'} ]\\
@V  \nabla_t VV @V \:\: \nabla_t VV @V  \nabla_t V \:\:V \\
\Omega^1_{\bC}\otimes \DR_{/\bC}({\mathcal M'}) @=
[ \Omega^1_{\bC}\otimes \Omega^k_{M'/\bC}\otimes {\mathcal M'} @> \nabla_{/\bC} >>   \Omega^1_{\bC} \otimes\Omega^{k+1}_{M'/\bC}\otimes  {\mathcal M'} ] \:,
 \end{CD}$$
with the ``top-dimensional forms" $\Omega^1_{\bC} \otimes\Omega^{m}_{M'/\bC}$ in bidegree
$(0,0)$. Here the horizontal lines come from the corresponding {\em relative de Rham complex} 
$\DR_{/\bC}({\mathcal M'})$ of 
${\mathcal M'}$ viewed only as left ${\mathcal D}_{M'/\bC}$-module, with ${\mathcal D}_{M'/\bC}\subset  V_{0}{\mathcal D}_{M'}$.

Moreover, $\DR({\mathcal M'})$  becomes a bifiltered double complex by
{\footnotesize
$$
\begin{CD} F_pV_{\alpha}\DR_{/\bC}({\mathcal M'}) @=
[ \Omega^k_{M'/\bC} \otimes F_{p+k}V_{\alpha}{\mathcal M'} @> \nabla_{/\bC} >>   \Omega^{k+1}_{M'/\bC} \otimes F_{p+k+1}V_{\alpha}{\mathcal M'} ]\\
@V  \nabla_t VV @V \:\: \nabla_t VV @V  \nabla_t V \:\:V \\
\Omega^1_{\bC}\otimes F_{p+1}V_{\alpha+1}\DR_{/\bC}({\mathcal M'}) @=
[  \Omega^1_{\bC}\otimes \Omega^k_{M'/\bC}\otimes  F_{p+k+1}V_{\alpha+1}{\mathcal M'} @> \nabla_{/\bC} >>   \Omega^1_{\bC} \otimes\Omega^{k+1}_{M'/\bC}\otimes   F_{p+k+2}V_{\alpha+1}{\mathcal M'} ] 
 \end{CD}
$$
}

Again, all differentials in the $F$-graded pieces of the $V$-filtered complex
\begin{equation}
\begin{CD} Gr^F_pV_{\alpha}\DR_{/\bC}({\mathcal M'})  \\
@V  Gr^F(\nabla_t) VV\\
\Omega^1_{\bC}\otimes Gr^F_{p+1}V_{\alpha+1}DR_{/\bC}({\mathcal M'})
  \end{CD}
\end{equation}
are ${\mathcal O}_{M'}$-linear.

From the strict specializability property (s3), one can show the following:
\bp 
The horizontal inclusion of $F$-filtered double complexes
$$\begin{CD} F_pV_{-1}\DR_{/\bC}({\mathcal M'}) @>>> F_p\DR_{/\bC}({\mathcal M'}) \\
@V  \nabla_t VV  @V  \nabla_t VV\\
\Omega^1_{\bC}\otimes F_{p+1}V_{0}\DR_{/\bC}({\mathcal M'}) @>>>
\Omega^1_{\bC}\otimes F_{p+1}\DR_{/\bC}({\mathcal M'})
  \end{CD}$$
induces a {\em filtered quasi-isomorphism} of the corresponding total complexes.
\ep

So the total complex of 
\begin{equation}
\begin{CD} Gr^F_pV_{-1}\DR_{/\bC}({\mathcal M'})  \\
@V  Gr^F_p(\nabla_t) VV \\
\Omega^1_{\bC}\otimes Gr^F_{p+1}V_{0}\DR_{/\bC}({\mathcal M'})
  \end{CD}
  \end{equation}
  represents $Gr^F_p \DR({\mathcal M'})$ and it is acyclic for almost all $p$. Hence it can be used for the calculation of 
  $$\DR_y([{\mathcal M'}]) \quad  \text{and} \quad i^! \DR_y([{\mathcal M'}]) \:.$$
Moreover, it belongs to $D^b_{{\rm coh}}(i'(M))$ because $\cM'$ is a mixed Hodge module supported on $i'(M)$, so that all $Gr^F_p\cM'$ are in fact $\cO_{i'(M)}$-modules by \cite{Sa0}[Lem.3.2.6].

For the calculation of  $i^! \DR_y([{\mathcal M'}])$ we can even use the total complex of
\begin{equation}
\begin{CD} Gr^F_pV_{<-1}\DR_{/\bC}({\mathcal M'})  \\
@V  Gr^F_p(\nabla_t) VV \\
\Omega^1_{\bC}\otimes Gr^F_{p+1}V_{<0}\DR_{/\bC}({\mathcal M'}) 
  \end{CD}
  \end{equation}
because 
$$
Gr^F_p Gr^V_{-1}\DR_{/\bC}({\mathcal M'}) 
\simeq i_*( Gr^F_p \DR(Gr^V_{-1}{\mathcal M'}|_M) ) 
$$
and 
$$
Gr^F_{p+1} Gr^V_{0}\DR_{/\bC}({\mathcal M'})  \simeq i_*(Gr^F_{p+1} \DR(Gr^V_{0}{\mathcal M'}|_M)).
$$
Indeed, as the right-hand complexes $Gr^F_p\DR$ belong to $D^b_{\rm coh}(M)$ by property (s1), they have well-defined classes in $K_0(M)$,  and moreover  $i^!i_*=0$ for the global hypersurface inclusion $i:M=\{t=0\} \to M'$.

Finally, we get:
\begin{equation}\begin{split}
Li^* Gr^F_{p} V_{<-1}\DR_{/\bC}( {\mathcal M'}) &\simeq 
i^* Gr^F_{p} V_{<-1} \DR_{/\bC}( {\mathcal M'})\\ 
&\overset{\cdot t}{\simeq} Gr^F_{p} \DR(\;V_{<0}/ tV_{<0}{\mathcal M'}|_M\;)\\
&\simeq Gr^F_{p} \DR(\;V_{<0}/ V_{<-1}{\mathcal M'}|_M\;),
\end{split}\end{equation}
where the first isomorphism uses the fact that there is no $t$-torsion (by (\ref{star}) and inductive use of property (s2)), and the last two isomorphisms follow from property (s2).
Similarly, we have:
\begin{equation}
Li^* Gr^F_{p+1} V_{<0}\DR_{/\bC}( {\mathcal M'}) \simeq 
Gr^F_{p+1} \DR(\;V_{<0}/ V_{<-1}{\mathcal M'}|_M\;).
\end{equation}

Putting everything together, we get  by additivity the following equality in $K_0(X)$:
\begin{equation}\begin{split}
&i^![Gr^F_p \DR({\mathcal M'})] =\\ -[Gr^F_{p} \DR(\;V_{<0}/ V_{<-1}&{\mathcal M'}|_M\;)]
+ [Gr^F_{p+1} \DR(\;V_{<0}/ V_{<-1}{\mathcal M'}|_M\;)]  \:.
\end{split}\end{equation}
Note that the minus sign for the first class on the right side comes from the fact
that $Gr^F_{p} \DR(\;V_{<0}/ V_{<-1}{\mathcal M'}|_M\;)$, when regarded as a subcomplex of the corresponding double complex, agrees only up to a shift by $1$ with the usual convention that ``top-dimesional form'' are in degree zero. 

By using the filtration $V_{\beta}/ V_{<-1}{\mathcal M'}$ of
$V_{<0}/ V_{<-1}{\mathcal M'}$ by ${\mathcal D}_{M'/\bC}$-modules
($-1\leq \beta <0$), we get by property (s1) and additivity
the following equality in the Grothendieck group $K_0(X)$:
\begin{equation}\begin{split}
&i^![Gr^F_p \DR({\mathcal M'})] =\\ \sum_{-1\leq \beta <0}\:(\:
-[Gr^F_{p} \DR(\;&Gr^V_{\beta}{\mathcal M'}|_M\;)]
+ [Gr^F_{p+1} \DR(\;Gr^V_{\beta}{\mathcal M'}|_M\;)]\:) \:.
\end{split}\end{equation}
This implies Theorem \ref{motthm}.


\subsection{Application: Hirzebruch-Milnor Classes of Singular Hypersurfaces}\label{applic}

Let $X=\{f=0\}$ be an algebraic variety defined as the zero-set of codimension one of an algebraic function $f:M \to \bC$, for $M$ a  complex algebraic manifold of complex dimension $n+1$.  Let $i:X \hookrightarrow M$ be the inclusion map. Denote by $L\vert_X$ the trivial line bundle on $X$. Then the virtual tangent bundle of $X$ can be identified with 
\begin{equation}\label{vir} T^{\rm vir}_X=[T_M\vert_X - L\vert_X]\:,\end{equation}
since $N_XM\simeq f^*N_{\{0\}}\bC\simeq L\vert_X$.

Let as before $$\psi_f^H, \varphi_f^H: \MHM(M) \to \MHM(X)$$ be the nearby and resp. vanishing cycle functors associated to $f$,  defined on the level of Saito's algebraic mixed Hodge modules. These functors induce transformations on the corresponding Grothendieck groups and, by construction,  the following identity holds in $K_0(\MHM(X))$ for any $[\cM]\in K_0(\MHM(M))$:
\begin{equation}\label{triangle}
\psi^H_f([\cM])= \varphi^H_f([\cM])-i^*([\cM])\: .
\end{equation}
Recall that the shifted transformations $\psi'^H_f:=\psi^H_f[1]$ and $\varphi'^H_f:=\varphi^H_f[1]$
correspond under the forgetful functor ${\rat}$ to the usual nearby and vanishing cycle functors.

Let $i^!:H_*(M) \to H_{*-1}(X)$ denote the Gysin map between the corresponding  homology theories (see \cite{Fu,V}). 
An easy consequence of the 
{specialization property} (\ref{sp}) for the Hirzebruch class transformation is the following:
\bl
\begin{equation}\label{eq10} {\widehat T}_{y*}^{\rm vir}(X):=\widehat{T}_y^*(T^{{\rm vir}}_X) \cap [X]=\widehat{T}_{y*}(\psi'^H_f(\left[\bQ^H_M\right])).
\end{equation}
\el
\begin{proof} Since $M$ is smooth, it follows that $\bQ^H_M$ is a shifted mixed Hodge module. 
By applying the identity (\ref{sp}) to the class $[\bQ_M^H]\in K_0(\MHM(M))$ we have that 
$$\widehat{T}_{y*}(\psi'^H_f(\left[\bQ^H_M\right]))=i^! \widehat{T}_{y*}([\bQ^H_M])=i^!\widehat{T}_{y*}(M)=i^!({\widehat T}_y^*(T_M) \cap [M])\:,$$
where the last identity follows from the normalization property of Hirzebruch classes as $M$ is smooth. Moreover, by the definition of the Gysin map, the last term of the above identity becomes 
$$ i^*(\widehat{T}_y^*(T_M)) \cap i^![M] = i^*(\widehat{T}_y^*(T_M)) \cap [X]\:,$$ 
which by the identification in (\ref{vir}) is simply equal to $\widehat{T}_{y*}^{\rm vir}(X)$.
\end{proof}

We thus have the following:
\bt\label{M} Let $X=f^{-1}(0)$ be a globally defined hypersurface (of codimension one) in a complex algebraic manifold $M$. Then the difference class 
$$\cM\widehat{T}_{y*}(X):=\widehat{T}_{y*}^{\rm vir}(X) - \widehat{T}_{y*}(X)$$ is entirely determined by the vanishing cycles of $f:M \to \bC$. 
More precisely,  
\begin{equation}\label{eq20}
\cM\widehat{T}_{y*}(X)=\widehat{T}_{y*}(\varphi'^H_f(\left[\bQ^H_M\right])).
\end{equation}
\et

Since the complex $\varphi^H_f(\bQ_M)$ is supported only on the singular locus $X_{\rm sing}$ of $X$ (i.e., on the set of points in $X$ where the differential $df$ vanishes), the result of Theorem \ref{M} shows that the difference class $\cM\widehat{T}_{y*}(X):=\widehat{T}_{y*}^{\rm vir}(X) - \widehat{T}_{y*}(X)$ can be expressed entirely only in terms of invariants of the singularities of $X$. Namely, by the functoriality of the transformation $\widehat{T}_{y*}$ (for the closed inclusion $X_{\rm sing}\hookrightarrow X$), we can view
$\cM\widehat{T}_{y*}(X) \in H_*(X_{\rm sing})\otimes\bQ[y]$ 
as a {\it localized} class.
Therefore, we have the following
\bc
The classes $\widehat{T}_{y*}^{\rm vir}(X)$ and $\widehat{T}_{y*}(X)$ coincide in dimensions greater than the dimension of the singular locus of $X$, i.e.,
$$\widehat{T}_{y,i}^{\rm vir}(X)=\widehat{T}_{y,i}(X) \in H_i(X)\otimes \bQ[y] \quad \text{for \ $i> \dim(X_{\rm sing})$ .}$$
\ec

By using the sets of generators of $K_0(\MHM(X))$, as described in Remark \ref{generator}, we can obtain precise formulae for the difference class $\cM\widehat{T}_{y*}(X):=\widehat{T}_{y*}^{\rm vir}(X) - \widehat{T}_{y*}(X)$ of a globally defined hypersurface. For simplicity, we assume here that the monodromy contributions along all strata in a stratification of $X$ are trivial, e.g., all strata are simply-connected (otherwise, one would have to use Proposition \ref{39}). This assumption, together with the rigidity property for variations and the multiplicativity for external products,  allow us to identify the coefficients in the above generating sets of $K_0(\MHM(X))$. For example, we have the following result:

\bt\label{main} Let $X=\{f=0\}$ be a complex algebraic variety defined as the zero-set (of codimension one)  of an algebraic function $f:M \to \bC$, for $M$ a  complex algebraic manifold.
Let $\cS_0$ be a partition of the singular locus $X_{\rm sing}$ into disjoint locally closed complex algebraic submanifolds $S$, such that the restrictions $\varphi_f(\bQ_M)\vert_S$ of the vanishing cycle complex to all pieces $S$ of this partition have constant cohomology sheaves (e.g., these are locally constant sheaves on each $S$, and the pieces $S$ are simply-connected).
For each $S \in \cS_0$, let $F_s$ be the Milnor fiber of a point $s \in S$. Then:
\begin{equation}\label{eq.main-part1}
\widehat{T}_{y*}^{\rm vir}(X) - \widehat{T}_{y*}(X)
=\sum_{S \in \cS_0} \;\left( \widehat{T}_{y*}(\bar S) - \widehat{T}_{y*}({\bar S} \setminus S) \right) \cdot \chi_y([\widetilde{H}^*(F_s;\bQ)]) \:.
\end{equation} 
\et

For example, if $X$ has only {\it isolated} singularities,
 the two classes $\widehat{T}_{y*}^{\rm vir}(X)$ and $\widehat{T}_{y*}(X)$ coincide except in degree zero, where their difference is measured (up to a sign) by the sum of Hodge polynomials associated to the middle cohomology of the corresponding Milnor fibers attached to the singular points. More precisely, we have in this case that:
\begin{equation}
\widehat{T}_{y*}^{\rm vir}(X) - \widehat{T}_{y*}(X)=\sum_{x \in X_{\rm sing}} (-1)^{n} \chi_y([\widetilde{H}^n(F_x;\bQ)])= \sum_{x\in X_{\rm sing}} \chi_y([\widetilde{H}^*(F_x;\bQ)]),
\end{equation}
where  $F_x$ is the Milnor fiber of the isolated hypersurface singularity germ $(X,x)$, and $n$ is the complex dimension of $X$. 
Recall that the cohomology groups $\widetilde{H}^k(F_s;\bQ)$  carry  canonical mixed Hodge structures 
 coming from the stalk formula
\begin{equation}\label{stalk}
\widetilde{H}^k(F_s;\bQ)\simeq H^k(\varphi_f(\bQ_M)_s).
\end{equation}
In the case of isolated singularities, the Hodge $\chi_y$-polynomials of the Milnor fibers at singular points are just Hodge-theoretic refinements of the Milnor numbers since
$$\chi_{-1}([\widetilde{H}^*(F_x;\bQ)])= \chi([\widetilde{H}^*(F_x;\bQ)])$$
is the reduced Euler characteristic of the Milnor fiber $F_x$, which up to a sign is just the Milnor number at $x$.
 For this reason, we regard the difference 
 \begin{equation}
 \cM\widehat{T}_{y*}(X):=\widehat{T}_{y*}^{\rm vir}(X) - \widehat{T}_{y*}(X) \in H_*(X) \otimes \bQ[y]
 \end{equation} 
 as a Hodge-theoretic Milnor class, and call it the \emph{Hirzebruch-Milnor class} of the hypersurface $X$. In fact, it is always the case that by substituting $y=-1$ into $\cM\widehat{T}_{y*}(X)$ we obtain the (rational) Milnor class ${\cM}_*(X) \otimes \bQ$ of $X$.

\br The above theorem can be used for computing the homology Hirzebruch class of the Pfaffian hypersurface and, respectively,  of the Hilbert scheme $$(\bC^3)^{[4]}:=\{df_4=0\}\subset M_4$$ considered in \cite{DS}[Sect.2.4 and Sect.3]. Indeed, the singular loci of the two hypersurfaces under discussion have ``adapted" partitions as in the above theorem, with only simply-connected strata (cf. \cite{DS}[Lem.2.4.1 and Cor.3.3.2]). Moreover, the mixed Hodge module corresponding to the vanishing cycles of the defining function, and its Hodge-Deligne polynomial are calculated in \cite{DS}[Thm.2.5.1, Thm.2.5.2, Cor.3.3.2 and Thm.3.4.1].
So Theorem \ref{main} above can be used for obtaining class versions of these results from \cite{DS}.
\er

\br For similar (inductive) calculations in more general situations, e.g., projective or global complete intersections, see \cite{MSS2,Sd}.
\er


\section{Equivariant Characteristic Classes}\label{equivariant}


\subsection{Motivation. Construction. Properties}

\subsubsection{Motivation} 

Equivariant characteristic numbers (genera) of complex algebraic varieties are generally defined by combining the information encoded by the filtrations of the mixed Hodge structure in cohomology with the action of a finite group preserving these filtrations (e.g., an algebraic action). 

For example, the equivariant Hirzebruch polynomial $\chi_y(X;g)$ considered here is defined in terms of the Hodge filtration on the (compactly supported) cohomology of a complex algebraic variety $X$ acted upon by a finite group $G$ of algebraic automorphisms $g$ of $X$. More precisely, we define
$$\chi_y(X;g):=\sum_{i,p} (-1)^i trace\left(g|Gr^p_FH_{({c})}^i(X;\bC)\right) \cdot (-y)^p.$$

One of the main motivations for studying such equivariant invariants is the information they provide when comparing invariants of an algebraic variety to those of its orbit space. For example, the equivariant Hirzebruch genera $\chi_y(X;g)$, $g \in G \setminus  \{id\}$, of a projective variety $X$ measure the ``difference" between the Hirzebruch polynomials $\chi_y(X)$ and, resp., $\chi_y(X/G)$. More precisely, in \cite{eCMS} we prove that:
$$\chi_y(X/G)=\frac{1}{|G|} \sum_{g \in G} \chi_y(X;g).$$
A similar relationship was used by Hirzebruch in order to compute the signature of certain ramified coverings of closed manifolds.
This calculation suggests the following {\it principle} which is obeyed by many invariants of global quotients:\newline {\it ``If $G$ is a {finite} group acting algebraically on a complex quasi-projective variety $X$, invariants of the orbit space $X/G$ are computed by an appropriate averaging of {equivariant invariants} of $X$"}.\\
So computing invariants of global quotients translates into the computation of equivariant invariants.

If $X$ is a compact algebraic {\it manifold}, the {\it Atiyah-Singer holomorphic Lefschetz formula} \cite{AS,HZ} can be used to compute the equivariant Hirzebruch polynomial $\chi_y(X;g)$ in terms of characteristic classes of the fixed point set $X^g$ and of its normal bundle in $X$: 
\be \chi_y(X;g)=\int_{[X^g]} T^*_y(X;g) \cap  [X^g]  , \ee with $T^*_y(X;g) \in H^*(X^g) \otimes \bC[y]$ the {\it cohomological Atiyah-Singer class}.

In the singular setup, in \cite{CMSS2} we define homological Atiyah-Singer classes $$T_{y*}(X;g) \in H_*(X^g) \otimes \bC[y]$$ for singular quasi-projective varieties, and prove a singular version of the Atiyah-Singer holomorphic Lefschetz formula for projective varieties:
\be \chi_y(X;g)=\int_{[X^g]} T_{y*}(X;g).\ee
(Here $H_*(-)$ denotes as before the Borel-Moore homology in even degrees.) The 
construction and its applications shall be explained below.


\subsubsection{Construction} Let $X$ be a (possibly singular) quasi-projective variety acted upon by a finite group $G$ of algebraic automorphisms.
The Atiyah--Singer class transformation $${T_y}_*(-; g):K_0(\MHM^G(X)) \to H_*(X^g) \otimes \bC[y^{\pm 1}]$$  is constructed in \cite{CMSS2} in two stages. First, by using  Saito's theory of algebraic mixed Hodge modules, we construct an equivariant version of the Hodge-Chern class transformation of \cite{BSY}, i.e., 
\begin{equation}
\DR_y^G:K_0(\MHM^G(X)) \to K^G_0(X) \otimes \bZ[y^{\pm 1}],
\end{equation}
for $K^G_0(X) := K_0({\rm Coh}^G(X))$ the Grothendieck group of $G$-equivariant algebraic coherent sheaves on $X$.
Secondly, we employ  the Lefschetz--Riemann--Roch transformation of Baum--Fulton--Quart \cite{BFQ} and Moonen \cite{Mo},
\begin{equation}
td_*(-; g):K^G_0(X)  \to H_*(X^g)\otimes\bC,
\end{equation}
to obtain (localized) homology classes on the fixed-point set $X^g$. 

\bigskip

In order to define the equivariant Hodge-Chern class transformation $\DR_y^G$, we work in the category $D^{b,G} \MHM(X)$ of $G$-equivariant objects in the derived category $D^b\MHM(X)$ of algebraic mixed Hodge modules on $X$, and similarly for $D_{\rm coh}^{b,G}(X)$,  the category of $G$-equivariant objects in the derived category $D_{\rm coh}^{b}(X)$ of bounded complexes of $\cO_X$-sheaves with coherent cohomology. Let us recall that in both cases,  a $G$-equivariant element $\cM$ is just an element in the underlying additive category (e.g., $D^b\MHM(X)$), with a $G$-action given by isomorphisms
$$\mu_g: \cM \to g_* \cM \quad (g\in G),$$  
such that $\mu_{id}=id$ and $\mu_{gh}=g_*(\mu_h)\circ \mu_g$ for all $g,h \in G$ (see 
\cite{MS09}[Appendix]). These ``weak equivariant derived categories'' $D^{b,G}(-)$ are not triangulated in general, hence they are different from the Berstein-Lunts notion of an equivariant derived category \cite{BL}. Nevertheless, one can define a suitable Grothendieck group, by using ``equivariant distinguished triangles'' in the underlying derived category $D^b(-)$, and obtain  isomorphisms (cf. \cite{CMSS2}[Lemma 6.7]): 
$$ K_0(D^{b,G} \MHM(X))=K_0(\MHM^G(X)) \quad \text{and} \quad K_0(D_{\rm coh}^{b,G}(X))=K^G_0(X).$$
The equivariant Hodge-Chern class transformation $\DR_y^G$ is defined by noting that Saito's natural transformations of triangulated categories (cf. Theorem \ref{grDR}) 
$$Gr^F_p \DR:D^b\MHM(X) \to D^b_{\rm coh}(X)$$ commute with the push-forward $g_*$ induced by each $g \in G$, thus inducing equivariant transformations (cf. \cite{CMSS2}[Example 6.6]) $$Gr^F_p \DR^G:D^{b,G}\MHM(X) \to D^{b,G}_{\rm coh}(X).$$ 
This simple definition of the equivariant transformation $Gr^F_p \DR^G$  depends crucially on our use of weak equivariant derived categories (and it is not clear how to obtain such a transformation in the singular setting by using the Bernstein-Lunts approach).
Note that for a fixed $\cM \in D^{b,G} \MHM(X)$, one has that $Gr^F_p \DR^G(\cM)=0$ for all but finitely many $p \in \bZ$.  This yields the following definition (where we use as before the notation $Gr^p_F$ in place of $Gr^F_{-p}$, corresponding to the identification $F^p=F_{-p}$):
\bd The {\it $G$-equivariant Hodge-Chern class transformation} 
$$\DR_y^G:K_0(\MHM^G(X)) \to K_0(D_{\rm coh}^{b,G}(X)) \otimes \bZ[y^{\pm 1}] = 
K^G_0(X) \otimes \bZ[y^{\pm 1}]$$
is defined by:
\begin{equation}\label{d}
\DR_y^G([\cM]):= \sum_{p} \left[Gr_F^{p} \DR^G (\cM) \right] \cdot (-y)^p=  \sum_{i,p}(-1)^i \left[\cH^i(Gr^p_F \DR^G (\cM))\right] \cdot (-y)^p .
\end{equation} 
The (un-normalized) {\it Atiyah--Singer class transformation} is defined by the composition
\begin{equation}
{T_y}_*(-; g):=td_*(-; g) \circ \DR_y^G,
\end{equation}
with 
\begin{equation}
td_*(-; g):K^G_0(X) \to H_{*}(X^g) \otimes \bC
\end{equation}
the Lefschetz--Riemann--Roch transformation (extended linearly over $ \bZ[y^{\pm 1}]$).
A normalized Atiyah--Singer class transformation can be defined similarly by using a twisted Lefschetz--Riemann--Roch transformation $td_{(1+y)*}(-; g)$ defined as before by multiplication by $(1+y)^{-k}$ on the degree $k$ component. The corresponding Atiyah-Singer characteristic classes ${T_y}_*(X; g)$ of a variety $X$ are obtained by evaluating the above transformation ${T_y}_*(-; g)$ at the class of the constant $G$-equivariant Hodge sheaf $\bQ^H_X$.
\ed

\subsubsection{Properties}

By construction, the transformations $\DR_y^G$ and ${T_y}_*(-; g)$ commute with proper pushforward and restriction to open subsets. Moreover, for a subgroup $H$ of $G$, these transformations commute with the obvious restriction functors ${\rm Res}^G_H$. For the trivial subgroup, this is just the forgetful functor ${\it For}:={\rm Res}^G_{\{id\}}$. If $G=\{id \}$ is the trivial group, then $\DR_y^G$ is just the (non-equivariant) Hodge-Chern class transformation of \cite{BSY,Sc}, while ${T_y}_*(-; id)$ is the complexified version of the (un-normalized) Hirzebruch class transformation ${T_y}_*$.

Our approach based on weak equivariant complexes of mixed Hodge modules allows us to formally extend most of the properties of Hodge-Chern and resp. Hirzebruch classes from \cite{BSY} and \cite{Sc}[Sect.4,5]  to the equivariant context. We give a brief account of these properties.
For example, we have:
\bp\label{db} Let $G$ be a finite group of algebraic automorphisms of a complex quasi-projective variety $X$, with at most Du Bois singularities.  Then \be \DR^G_0([\bQ^H_X])=[\cO_X] \in K^G_0(X),\ee as given by the class of the structure sheaf with its canonical $G$-action. In particular, \be{T_0}_*(X; g)=td_*([\cO_X];g)=:td_*(X;g).\ee
\ep
Similarly, it can be shown that both transformations are multiplicative with respect to external products (see \cite{CMSS2}[Thm.3.6,Cor.4.3]).

By the Lefschetz-Riemann-Roch theorem of \cite{BFQ,Mo}, the class transformation  
$td_*(-;g)$ commutes with the pushforward under proper morphisms, 
so the same is true for  the Atiyah-Singer transformations $T_{y*}(-;g)$ and $\widehat{T}_{y*}(-;g)$. If we apply this observation to the constant map $k:X\to pt$ with $X$ compact, then the
pushforward for $H_*$ is identified with the degree map, and we have that 
$K^G_0(pt)$ is the complex representation ring of $G$, and $\MHM^G(pt)$ is identified with the abelian category of $G$-equivariant graded-polarizable mixed $\bQ$-Hodge structures.
By definition, we have for $H^{\bullet}\in D^{b,G}\MHS^p$ that
\be T_{y*}(H^{\bullet};g)=\widehat{T}_{y*}(H^{\bullet};g)=\chi_y(H^{\bullet};g):=
\sum_{j,p}\,(-1)^j trace(g|\Gr_F^pH_{\bC}^j)\,(-y)^p.\ee

We also mention here the relation between the equivariant and resp. non-equivariant Hodge-Chern class transformations for spaces with trivial $G$-action. 
Let $G$ act trivially on the quasi-projective variety $X$. Then one can 
consider the projector $$(-)^G:=\frac{1}{\vert G \vert} \sum_{g \in G} \mu_g$$ acting on the categories $D^{b,G}\MHM(X)$ and $D^{b,G}_{\rm coh}(X)$, for $\mu_g$ the isomorphism induced from the action of $g \in G$. Here we use the fact that the underlying categories $D^{b}\MHM(X)$ and $D^{b}_{\rm coh}(X)$ are $\bQ$-linear additive categories which are moreover Karoubian (i.e., any projector has a kernel, see \cite{MS09} and the references therein).  
Since $(-)^G$ is exact, we obtain induced functors on the Grothendieck groups:
$$[-]^G:K_0 (\MHM^G(X)) \to K_0 (\MHM(X))$$ and
$$[-]^G:K^G_0(X) \to K_0(X).$$
We then have the following result:
\bp\label{L1}
Let $X$ be a complex quasi-projective $G$-variety, with a trivial action of the finite group $G$. Then the following diagram commutes:
\begin{equation}\label{LL1}
\begin{CD}
K_0(\MHM^G(X)) @> \DR_y^G >> K_0^G(X) \otimes \bZ[y^{\pm 1}]\\
@V [-]^G VV @VV [-]^G  V\\
K_0(\MHM(X)) @> \DR_y >> K_0(X) \otimes \bZ[y^{\pm 1}]\
\end{CD}
\end{equation}
\ep

This relationship can be used for computing characteristic classes of global quotients. Indeed, let $X':=X/G$ be the quotient space, with finite projection map $\pi:X \to X'$, which is viewed as a $G$-map with trivial action on $X'$. Then, it follows from \cite{CMSS2}[Lem.5.3] that:
\be\label{100}
\left(\pi_* \pi^*\cM \right)^G \simeq \cM,
\ee for any $\cM \in D^b\MHM(X')$, 
where $\pi^*\cM$ and resp. $\pi_* \pi^*\cM$ are considered with their respective induced $G$-actions.
Together with Proposition \ref{L1}, this yields the following result:
\bt\label{gen} Let $G$ be a finite group acting by algebraic automorphisms on the complex quasi-projective variety $X$.  Let  $\pi^g:X^g \to X/G$ be the composition of the projection map $\pi:X \to X/G$ with the inclusion $i^g:X^g \hookrightarrow X$. Then, for any $\cM \in D^{b,G}\MHM(X)$, we have:
\be\label{arbn}
{T_y}_*([\pi_*\cM]^G)=\frac{1}{\vert G \vert} \sum_{g \in G} \pi^g_* {T_y}_*([\cM];g).
\ee
\et
In particular, for $\cM=\bQ_X^H$, we get the following consequence:
\bc\label{orb} Let $G$ be a finite group acting by algebraic automorphisms on the complex quasi-projective variety $X$.  Let  $\pi^g:X^g \to X/G$ be the composition of the projection map $\pi:X \to X/G$ with the inclusion $i^g:X^g \hookrightarrow X$. Then
\begin{equation}\label{Z}
 {{T_y}}_*(X/G)=\frac{1}{\vert G \vert} \sum_{g \in G} \pi^g_*{{T_y}}_*(X;g).
\end{equation}
\ec

\br\label{Za}
Similar formulas hold for the normalized Atiyah-Singer classes. In particular, if $X'=X/G$ is the  quotient of a projective manifold by a finite algebraic group action, one gets by comparing the normalized version of formula (\ref{Z}) for $y=1$ with a similar $L$-class formula due to Moonen-Zagier \cite{Mo, Za}, the conjectured identity: $\widehat{T}_{1*}(X')=L_*(X')$. 
For more applications, e.g., defect formulas of Atiyah-Meyer type, see \cite{CMSS2}.
\er


\subsection{Symmetric Products of Mixed Hodge Modules}\label{symmhm}
Let $X$ be a (possibly singular) complex quasi-projective variety, and define its $n$-th {\it symmetric product} $X^{(n)}:=X^n/\Sigma_n$ as the quotient of the product of $n$ copies of $X$ by the natural action of the symmetric group of $n$ elements, $\Sigma_n$. The symmetric products of $X$ are quasi-projective varieties as well. Let $\pi_n:X^n \to \Xs$ be the natural projection map.

In \cite{MSS}, we define an action of the symmetric group $\Sigma_n$ on the $n$-fold external self-product $\boxtimes^n\cM$ of an arbitrary bounded complex of mixed Hodge modules $\cM \in D^b\MHM(X)$. By construction, this action is compatible with the natural action on the underlying $\bQ$-complexes. There are, however,  certain technical difficulties associated with this construction, since the difference in the $t$-structures
of the underlying $\cD$-modules and $\bQ$-complexes gives certain differences of signs. In \cite{MSS}[Prop.1.5,Thm.1.9] we solve this problem by showing that there is a sign cancellation. We also prove in \cite{MSS}[(1.12)] the equivariant K\"unneth formula for the $n$-fold external products of bounded complexes of mixed Hodge modules in a compatible way with the action of the symmetric group $\Sigma_n$.

\bd For a complex of mixed Hodge modules $\cM \in D^b\MHM(X)$, we  define its {\it $n$-th symmetric product} by:
\begin{equation}\label{e}{\cM^{(n)}}:=({\pi_n}_*\cM^{\boxtimes{n}})^{\Sigma_n} 
\in D^b\MHM(X^{(n)}),\end{equation} where  $\cM^{\boxtimes n} \in D^b\MHM(X^n)$ is the $n$-th external product of $\cM$ with the above mentioned $\Sigma_n$-action, and  $(-)^{\Sigma_n}$ is the projector on the $\Sigma_n$-invariant sub-object (which is well-defined since $D^b\MHM(\Xs)$ is a Karoubian $\bQ$-linear additive category). 
\ed

As special cases of (\ref{e}), it was shown in \cite{MSS}[Rem.2.4(i)] that for $\cM=\bQ_X^H$ the constant Hodge sheaf on $X$, one obtains:
\begin{equation}\label{q}\left(\bQ_X^H\right)^{(n)}=\bQ^H_{X^{(n)}}.\end{equation}
Moreover, it follows from the equivariant K\"unneth formula of \cite{MSS}[(1.12)] that for $\cM:=f_{*(!)}\bQ^H_Y$ with $f:Y \to X$ an algebraic map, we get 
\begin{equation}\label{qq}\left(f_{*(!)}\bQ_Y^H\right)^{(n)}=f^{(n)}_{*(!)}(\bQ^H_{Y^{(n)}}).\end{equation}

As a consequence of the above considerations, we obtain the following result (see \cite{MSS}[Thm.1]):
\bt
There is a canonical isomorphisms of graded mixed Hodge structures
$$H_{({c})}^{\bullet}(\Xs;{\cM^{(n)}})\simeq H_{({c})}^{\bullet}(X^n;\boxtimes^n\cM)^{\Sigma_n}
\simeq \bigl(\bigotimes^n H_{({c})}^{\bullet}(X;\cM)\bigr){}^{\Sigma_n},$$
in a compatible way with the corresponding isomorphisms
of the underlying $\bQ$-complexes.
\et
This result can be applied to prove a generating series formula for the mixed Hodge numbers of symmetric products of a given mixed Hodge module complex. More precisely, 
the Hodge numbers of $\cM\in D^b\MHM(X)$ are defined
for $p,q,k\in\bZ$ by
$$h_{({c})}^{p,q,k}(\cM):=h^{p,q}(H_{({c})}^k(X;\cM)):=
\dim_{\bC}(\Gr^p_F\Gr^W_{p+q}H_{({c})}^k(X;\cM)_{\bC}),$$
where $H_{({c})}^k(X;\cM)_{\bC}$ denotes the underlying $\bC$-vector space
of the mixed Hodge structure on $H_{({c})}^k(X;\cM)$.
Taking the alternating sums over $k$, we get the $E$-polynomial
in $\bZ[y^{\pm 1},x^{\pm 1}]$:
$$e_{({c})}(\cM;y,x):=\sum_{p,q} e_{({c})}^{p,q}(\cM)\,y^px^q,$$
where $$e_{({c})}^{p,q}(\cM):=\sum_k (-1)^k h_{({c})}^{p,q,k}(\cM).$$
Then the following result holds for the generating series of the above numbers and polynomials (see \cite{MSS}[Cor.2] and \cite{MS09}[Thm.1.1]):

\bt\label{57} For any bounded complex of mixed Hodge modules $\cM$ on a complex quasi-projective variety $X$, we have the following identities:
$$\aligned &\sum_{n \geq 0} \Bigl(\sum_{p,q,k} h_{({c})}^{p,q,k}({\cM^{(n)}})\,
y^px^q(-z)^k\Bigr)\,t^n=\prod_{p,q,k}\Bigl(\frac{1}{1-y^px^qz^kt}
\Bigr)^{(-1)^k h_{({c})}^{p,q,k}(\cM)},\\
&\sum_{n\geq 0}e_{({c})}({\cM^{(n)}};y,x)\,t^n=\prod_{p,q}\Bigl(
\frac{1}{1-y^px^qt}
\Bigr)^{ e_{({c})}^{p,q}(\cM)}= \exp\Bigl(\sum_{r \geq 1}e_{({c})}(\cM;y^r,x^r)\,
\frac{t^r}{r} \Bigr).\endaligned$$
\et
\br
The latter result, which holds for complex quasi-projective varieties with arbitrary singularities, includes many of the classical results in the literature as special cases. It also specializes to new generating series formulae, e.g., for the intersection cohomology Hodge numbers and resp. the Goresky-MacPherson intersection cohomology signatures of symmetric products of complex projective varieties, \cite{MSS}. For a different approach to more general generating series formulae, based on {\it pre-lambda rings}, see \cite{MS09}.
\er


\subsection{Application: Characteristic Classes of Symmetric Products}\label{applsym}
The standard approach for computing invariants of symmetric products of varieties (resp. mixed Hodge module complexes) is to encode the respective invariants of all symmetric products in a generating series, and to compute such an expression solely in terms of invariants of the original variety (resp. mixed Hodge module complex), e.g., see Theorem \ref{57} above.

To obtain generating series formulae for characteristic classes of symmetric products, one has to work in the corresponding {\it Pontrjagin homology ring}, e.g., see \cite{Mo} for Todd classes and \cite{Za} for $L$-classes.

In \cite{CMSSY}, we obtain the following generating series formula for the Hirzebruch classes of symmetric products of a mixed Hodge module complex (recall that $H_*$ denotes the Borel-Moore homology in even degrees):
\bt\label{genHir}
Let $X$ be a complex quasi-projective variety with associated Pontrjagin ring $$PH_*(X):=\sum_{n=0}^{\infty} \big( H_*(X^{(n)}) \otimes \bQ[y,y^{-1}] \big) \cdot t^n:=\prod_{n=0}^{\infty} H_*(X^{(n)}) \otimes \bQ[y,y^{-1}].$$
For any $\cM \in D^b\MHM(X)$, the following identity holds in $PH_*(X)$:
\begin{equation}\label{te1}
\sum_{n \geq 0} {{T_{(-y)}}_*} (\cM^{(n)} ) \cdot t^n= \exp \left( \sum_{r \geq 1} \Psi_r  \Big(d^r_* {{T_{(-y)}}_*}(\cM) \Big) \cdot \frac{t^r}{r} \right),\end{equation}
where
\begin{enumerate}
\item[(a)] $d^r:X \to X^{(r)}$ 
is the composition of  the diagonal embedding $i_r:X\simeq \Delta_r(X) \hookrightarrow X^r$ with the projection $\pi_r: X^r \to X^{(r)}$. 
\item[(b)] $\Psi_r$ is the \emph{$r$-th homological Adams operation}, which on $H_k(X^r) \otimes \bQ:=H^{BM}_{2k}(X^{{r}};\bQ)$ ($k \in \bZ$) is defined by multiplication by $\frac{1}{r^k}$, together with $y \mapsto y^r$.
\end{enumerate}
\et

Note that the diagonal maps $d^r$ in the above formula are needed to move homology classes from $X$ to the symmetric products $X^{({r})}$. The appearance on the homological Adams operation $\Psi_r$  will be explained later on.

In particular, for $\cM=\bQ_X^H$, formula (\ref{te1}) yields by (\ref{q}) a generating series formula for the Hirzebruch classes of symmetric products of a variety $X$:
\bc Let $X$ be a complex quasi-projective variety. Then in the above notations we have:
\begin{equation}\label{te100}
\sum_{n \geq 0} {{T_{(-y)}}_*} (X^{(n)}) \cdot t^n= \exp \left( \sum_{r \geq 1} \Psi_r  \left(d^r_* {{T_{(-y)}}_*}(X) \right) \cdot \frac{t^r}{r} \right).\end{equation}
\ec

\br
If $X$ is smooth and projective, formula (\ref{te100}) specializes to Moonen's generating series formula for his generalized Todd classes $\tau_y(X^{(n)})$ (cf. \cite{Mo}[p.172]). Moreover, after a suitable re-normalization (see \cite{CMSSY}[Sect.5]), formula (\ref{te100}) specializes for the value $y=1$ of the parameter to Ohmoto's generating series formula \cite{Oh} for the rational MacPherson--Chern classes $c_*(X^{(n)})$ of the symmetric products of $X$:
\begin{equation}\label{Oh}
\sum_{n \geq 0} {c_*} (X^{(n)} ) \cdot t^n= \exp \left( \sum_{r \geq 1}  d^r_* {c_*}(X)  \cdot \frac{t^r}{r} \right) .
\end{equation}
Formula (\ref{Oh}) is a characteristic class version of Macdonald's generating series formula for the topological Euler characteristic. More generally, formula (\ref{te1}) is a characteristic class version of the generating series formula for the Hodge polynomial $\chi_{-y}(-)=e(-;y,1)$ of Theorem \ref{57}.
\er

The strategy of proof of Theorem \ref{genHir} follows Zagier \cite{Za}[Ch.II] and Moonen \cite{Mo}[Ch.II, Sect.2], which deal with $L$-classes of symmetric products of a rational homology manifold and, respectively,  Todd classes of symmetric products of a complex projective variety. We give here a brief account of the steps involved in the proof.

Let $\sigma \in \Sigma_n$ have cycle partition $\lambda=(k_1, k_2, \cdots )$, i.e.,  $k_r$ is the number of length $r$ cycles in $\sigma$ and $n=\sum_r r \cdot k_r$. Let $$\pi^{\sigma}:(X^n)^{\sigma} \to X^{(n)}$$ denote the composition of the inclusion of the fixed point set $(X^n)^{\sigma} \hookrightarrow X^n$ followed by the projection $\pi_n:X^n \to X^{(n)}$. Then $\pi^{\sigma}$ is the product of projections $$d^r:X \simeq \Delta_r(X) \overset{i_r}{\hookrightarrow} X^r \overset{\pi_r}{\to} X^{(r)},$$ where each $r$-cycle contributes a copy of $d^r$. (Here $\Delta_r(X)$ denotes the diagonal in $X^r$.) 

Theorem \ref{genHir} is a consequence of the following sequence of reductions. 

First, by Theorem \ref{gen}, we get the following {\it averaging property}:

\bl\label{l1} For  $\cM \in D^b\MHM(X)$ and every $n \geq 0$, we have:
\begin{equation}\label{le1}
{T_{y}}_*(\cM^{(n)}) =\frac{1}{n!}\sum_{\sigma \in \Sigma_n} \pi^{\sigma}_*{T_{y}}_*(\cM^{\boxtimes{n}}; \sigma).
\end{equation}
\el
Secondly, the multiplicativity property of the Atiyah--Singer 
 class  transformation reduces the problem to the computation of the Atiyah-Singer classes corresponding to cycles:
\bl\label{l2} If $\sigma \in \Sigma_n$ has cycle-type $(k_1, k_2, \cdots)$, then:
\begin{equation}\label{le2}
{T_{y}}_*(\cM^{\boxtimes{n}}; \sigma)=\prod_r \left( {T_{y}}_*(\cM^{\boxtimes{r}}; \sigma_r) \right)^{k_r}.
\end{equation}
\el
Finally, the following {\it localization formula} holds:
\bl\label{l3} The following identification holds in $H_*(X) \otimes \bQ[y^{\pm 1}] \subset H_*(X) \otimes \bC[y^{\pm 1}]$:
\begin{equation}\label{le4}
{T_{(-y)}}_*(\cM^{\boxtimes{r}}; \sigma_r)=\Psi_r {T_{(-y)}}_*(\cM),
\end{equation}
with $\Psi_r$ the $r$-th homological Adams operation.
\el
This localization property is the key technical point in the proof of Theorem \ref{genHir}.
One ingredient in the proof of Lemma \ref{l3} is the corresponding Todd class formula due to Moonen \cite{Mo}[Satz 2.4, p.162]:
\bl\label{l33a}
Let $\sigma_r$  be an $r$-cycle. Then for any $\cG \in D^b_{\rm coh}(X)$, the following identity holds in $H_*(X) \otimes \bQ$:
\begin{equation}\label{e33a}
td_*(\cG^{\boxtimes r}; \sigma_r)=\Psi_r td_*(\cG)
\end{equation}
under the identification  $(X^r)^{\sigma_r} \simeq X$.
\el

\noindent The idea of proof uses an embedding $i: X \hookrightarrow M$ into a smooth complex algebraic variety $M$, together with a bounded locally free resolution $\cF$ of $i_*\cG$. So $i^r:X^r \hookrightarrow M^r$ is a $\Sigma_r$-equivariant embedding, with $\cF^{\boxtimes r}$ a $\Sigma_r$-equivariant locally free resolution of $(i_*\cG)^{\boxtimes r} \simeq i^r_*(\cG^{\boxtimes r})$. Then $\Delta_r^*(\cF^{\boxtimes r})$ can  be used for the calculation of $td_*(\cG^{\boxtimes r}; \sigma_r)$ in terms of a (suitably modified) localized Chern character, with $\Delta_r:M \to M^r$ the diagonal embedding. 
The homological Adams operation $\Psi_r$ of Lemma \ref{l33a} is induced from the $K$-theoretic Adams operation 
$$\Psi^r: K^0(M,M\setminus X) \to K^0(M,M\setminus X) \otimes \bC$$ appearing in Moonen's proof of his localization formula \cite{Mo}[p.164] as 
$$\Psi^r([\cF]):=[\Delta_r^*(\cF^{\boxtimes r})](\sigma_r).$$
Note that the $\langle\sigma_r\rangle$-equivariant vector bundle complex $\Delta_r^*(\cF^{\boxtimes r})$
is exact off $X$, i.e., it defines a class $$[\Delta_r^*(\cF^{\boxtimes r})] \in K_{
\langle\sigma_r\rangle}^0(M,M \setminus X) \simeq K^0(M,M \setminus X) \otimes R({
\langle\sigma_r\rangle}),$$ where $R({
\langle\sigma_r\rangle})$ is the complex representation ring of the cyclic group generated by $\sigma_r$. Then $$(\sigma_r):K_{
\langle\sigma_r\rangle}^0(M,M \setminus X) \simeq K^0(M,M \setminus X) \otimes R({
\langle\sigma_r\rangle}) \to K^0(M,M \setminus X) \otimes \bC$$
is induced by taking the trace homomorphism $tr(-;\sigma_r):R({
\langle\sigma_r\rangle})  \to \bC$ (see \cite{Mo}[p.67]).

The second ingredient needed for the proof of the localization formula (\ref{le4})  relies on understanding how Saito's functors $Gr^F_p\DR$ behave with respect to external products, so that Moonen's calculation can be adapted to the graded complexes $Gr^F_*\DR^{\Sigma_r}(\cM^{\boxtimes r})$, once we let $\Psi_r$ also act by $y \mapsto y^r$, see \cite{CMSSY}[Sect.4]. This is also the reason why we work with the classes ${T_{(-y)}}_*$.

\br The generating series formula for the Hirzebruch classes of symmetric products is further used in \cite{CMOSY} for proving a generating series formula for (the push-forward under the Hilbert-Chow morphism of) the Hirzebruch  classes of the {\it Hilbert schemes} of points for a quasi-projective manifold of arbitrary pure dimension.

The strategy outlined above for the generating series of Hirzebruch classes, is also used in \cite{CMSSY} for the study of Todd (resp. Chern) classes of symmetric products of coherent (resp. constructible) sheaf complexes, generalizing to arbitrary coefficients results of Moonen \cite{Mo} and resp. Ohmoto \cite{Oh}. 
\er

\end{document}